\documentclass{article}
\pdfoutput=1
\usepackage{amsmath}
\usepackage{amsfonts}
\usepackage{amssymb}
\usepackage{graphicx}
\usepackage{theorem}%
\newtheorem{theorem}{Theorem}

\newtheorem{corollary}[theorem]{Corollary}

\newtheorem{definition}[theorem]{Definition}

\newtheorem{proposition}[theorem]{Proposition}
\newtheorem{remark}[theorem]{Remark}

\theorembodyfont{\upshape}
\newtheorem{example}{Example}
\newenvironment{proof}[1][Proof]{\noindent\textbf{#1.} }{\ \rule{0.5em}{0.5em}}
\newenvironment{acknowledgement}[1][Acknowledgements]{\noindent\textbf{#1.} }{}
\begin{document}

\title{Riemannian groupoids and solitons for three-dimensional homogeneous Ricci and
cross curvature flows}
\author{David Glickenstein}
\maketitle
\tableofcontents

\section{Introduction}

In recent years, there has been significant progress towards understanding
geometric flows of Riemannian metrics, most notably with Hamilton's and,
later, Perelman's work on Ricci flow (See, e.g., \cite{H1}, \cite{H2},
\cite{P1}, \cite{P2}, \cite{KL}, \cite{CZ}, \cite{MT}). An important method in
the analysis of the Ricci flow is a careful classification and analysis of
singularities of solutions to the flow \cite{H2}, often using a compactness
theorem about manifolds with some kind of curvature bound (e.g., \cite{H3},
\cite{Lu1}, \cite{P1}, \cite{G1}, \cite{Lot1}). By using a compactness
theorem, one may extract a limit of solutions as the flow approaches a
singularity, and the limit gives information about the asymptotic behavior of
the flow.

The limit of a parabolic flow is expected to be highly symmetric, usually some
type of self-similar solution, also called a soliton. A soliton is a
generalization of a fixed point of a flow; in fact, it is a solution which is
a fixed point except that the metric could be changing by time dependent
diffeomorphisms and rescaling.

In this paper, we study three-dimensional homogeneous geometries. These
geometries are easier to study because we are able to describe the metric
explicitly and also exhibit a large number of diffeomorphisms to find limit
soliton metrics. In general, exhibiting these diffeomorphisms is likely to be
much more difficult. Three-dimensional homogeneous solutions of Ricci flow
were first studied by Isenberg and Jackson \cite{IJ} and later by Knopf-McLeod
\cite{KM}. The solutions of the simply connected homogeneous solutions were
described in some detail. These solutions are quite interesting since they
mostly exhibit a particular singularity type (Type III) and are often
collapsing with bounded curvature. Later, Lott \cite{Lot1} was able to use the
formalism of Riemannian groupoids to better understand the case of compact
homogenous geometries and gave a complete classification in dimension 3.

The purpose of this paper is to apply the techniques of Riemannian groupoids
to study the long term behavior of solutions of the negative cross curvature
flow (XCF), a geometric flow on three-manifolds first introduced by Chow and
Hamilton \cite{CH}. The behavior of the simply connected geometries was first
given by Cao, Ni, and Saloff-Coste \cite{CNS}. We explain what happens to
compact quotients of homogeneous solutions to XCF in a way similar to Lott's
work on Ricci flow. We include detailed analysis of the Ricci flow situation
as well, both to better explain our coordinates, a few of which differ from
Lott's treatment, and to emphasize the similarities in the techniques and
utility of the Riemannian groupoid formalism.

The rest of the paper is organized as follows. First we review the notions of
Riemannian groupoids. We then review relevant aspects of the Ricci flow and
cross curvature flow, together with theory of singularities and soliton
solutions. We then give detailed descriptions of the homogeneous geometries
$\operatorname{Nil},$ $\operatorname{Sol},$ $\widetilde{\operatorname{SL}%
}\left(  2,\mathbb{R}\right)  ,$ and $\widetilde{\operatorname{Isom}}\left(
\mathbb{E}^{2}\right)  .$ Note that we choose not to include the other
homogeneous geometries since they lack the complexity of these; that is, there
is no need to consider changing diffeomorphisms in their convergence.

\begin{acknowledgement}
The author would like to thank John Lott for a number of helpful conversations
about his work and for clarifying several of the constructions related to
Riemannian groupoids. The author would like to thank Xiaodong Cao for his help
with understanding the evolution of homogeneous geometries under XCF. The
author would also like to thank B. Chow, P. Foth, C. Guenther, J. Isenberg, D.
Knopf, and T. Payne for stimulating and helpful conversations on topics
related to this paper.
\end{acknowledgement}

\section{Riemannian Groupoids}

Haefliger first introduced the notion of Riemannian groupoid \cite{BH}. We
will primarily follow the exposition in \cite{Lot1}. In order to emphasize
pieces, we will re-introduce the definitions. We make an effort to provide the
minimal number of definitions to understand the statement of convergence. A
Riemannian groupoid is a structure that encapsulates the notions of a
manifold, orbifold, and quotient manifold in the same global definition. Two
excellent references for smooth groupoids are the books \cite{MM} and
\cite{Mac} (both sources refer to the smooth groupoids used here as Lie
groupoids). Riemannian groupoids were previously introduced in \cite{GGHR},
but the definition is slightly different (in Lott's treatment, one only needs
a Riemannian metric on $G^{\left(  0\right)  }$ and not on $G^{\left(
1\right)  }$).

Before we recall the definition, let us give two examples which give the
flavor of some things that a groupoid can do.

\begin{example}
[Manifold with charts]\label{manifold example}Let $\left\{  U_{i}\right\}
_{i\in I}$ be an open covering of a Riemannian manifold $M.$ The groupoid
perspective represents $M$ as two pieces, $G^{\left(  0\right)  }%
=\coprod\limits_{i\in I}U_{i}$, where $\coprod$ denotes the disjoint union,
and $G^{\left(  1\right)  },$ which consists of maps between the two points in
the disjoint union which correspond to the same point in the covering (i.e.,
if $x\in U_{i}\cap U_{j},$ there is a map $\left(  x_{i}\rightarrow
x_{j}\right)  \in G^{\left(  1\right)  }$ mapping the corresponding points in
the disjoint union). Note that there is always the identity map which maps a
point to itself, which we may consider as an embedding $e:G^{\left(  0\right)
}\rightarrow G^{\left(  1\right)  }.$ Every element in $G^{\left(  1\right)
}$ looks like $\left(  x\rightarrow y\right)  $ where $x,y\in G^{\left(
0\right)  },$ so there are source and range maps $s:G^{\left(  1\right)
}\rightarrow G^{\left(  0\right)  },$ $r:G^{\left(  1\right)  }\rightarrow
G^{\left(  0\right)  }$ that look like $s\left(  x\rightarrow y\right)  =x,$
$r\left(  x\rightarrow y\right)  =y.$ Furthermore, if the source of
$\gamma_{1}\in G^{\left(  1\right)  }$ is equal to the range of $\gamma_{2}\in
G^{\left(  1\right)  },$ e.g., $\gamma_{1}=\left(  y\rightarrow z\right)  $
and $\gamma_{2}=\left(  x\rightarrow y\right)  ,$ then there is a product
$\gamma_{1}\gamma_{2}$ which essentially is associativity, e.g.
\begin{align*}
\gamma_{1}\gamma_{2}  &  =\left(  y\rightarrow z\right)  \left(  x\rightarrow
y\right)  =x\rightarrow y\rightarrow z\\
&  =x\rightarrow z.
\end{align*}
Note that this only works because the $y$ is the same point, otherwise this
product is not defined. Furthermore, there are inverses, $\left(  x\rightarrow
y\right)  ^{-1}=\left(  y\rightarrow x\right)  .$
\end{example}

\begin{example}
[Quotient by a group action]\label{quotient example}Let $\Gamma$ be a group
acting on a space $X$ from the right. We will consider a groupoid structure
that represents the quotient space $X/\Gamma.$ Here we let $G^{\left(
0\right)  }=X$ and $G^{\left(  1\right)  }=\bigcup\limits_{\gamma\in
\Gamma,x\in G^{\left(  0\right)  }}\left\{  \left(  x\rightarrow
x\gamma\right)  \right\}  .$ It is easy to see that the maps $e,$ $s,$ $r$ are
defined here in much the same way as Example \ref{manifold example}. The group
action hypotheses imply that the product is well defined.
\end{example}

Let us recall the general definition of a groupoid:

\begin{definition}
\label{groupoid def}A groupoid $G=\left(  G^{\left(  0\right)  },G^{\left(
1\right)  },e,s,r,\cdot\right)  $ is a 6-tuple such that

\begin{enumerate}
\item $G^{\left(  0\right)  }$ and $G^{\left(  1\right)  }$ are sets.

\item The unit map $e:G^{\left(  0\right)  }\rightarrow G^{\left(  1\right)
}$ is an injection.

\item The source and range maps $s,r:G^{\left(  1\right)  }\rightarrow
G^{\left(  0\right)  }$ satisfy $s\circ e=r\circ e$ are the identity map.

\item The partially defined multiplication $\cdot:G^{\left(  1\right)  }\times
G^{\left(  1\right)  }\rightarrow G^{\left(  1\right)  },$ usually denoted by
juxtaposition, satisfies the following:

\begin{enumerate}
\item If $\gamma,\gamma^{\prime}\in G^{\left(  1\right)  },$ then the product
$\gamma\gamma^{\prime}$ is defined only if $s\left(  \gamma\right)  =r\left(
\gamma^{\prime}\right)  ;$ in this case, $s\left(  \gamma\gamma^{\prime
}\right)  =s\left(  \gamma^{\prime}\right)  $ and $r\left(  \gamma
\gamma^{\prime}\right)  =r\left(  \gamma\right)  .$

\item The product is associative, i.e., $\left(  \gamma\gamma^{\prime}\right)
\gamma^{\prime\prime}=\gamma\left(  \gamma^{\prime}\gamma^{\prime\prime
}\right)  ,$ if both sides make sense.

\item $\gamma e\left(  s\left(  \gamma\right)  \right)  =e\left(  r\left(
\gamma\right)  \right)  \gamma=\gamma.$

\item For any $\gamma\in G^{\left(  1\right)  },$ there is an element
$\gamma^{-1}\in G^{\left(  1\right)  }$ such that $s\left(  \gamma
^{-1}\right)  =r\left(  \gamma\right)  $, $r\left(  \gamma^{-1}\right)
=s\left(  \gamma\right)  ,$ $\gamma\gamma^{-1}=e\left(  r\left(
\gamma\right)  \right)  ,$ and $\gamma^{-1}\gamma=e\left(  s\left(
\gamma\right)  \right)  .$
\end{enumerate}
\end{enumerate}
\end{definition}

\begin{remark}
Elements of $G^{\left(  0\right)  }$ are called \emph{objects} and elements of
$G^{\left(  1\right)  }$ are called \emph{arrows}.
\end{remark}

\begin{definition}
A \emph{trivial groupoid} is one in which $G^{\left(  1\right)  }=G^{\left(
0\right)  }$ and $s$ and $r$ are both the identity map.
\end{definition}

Note that in Definition \ref{groupoid def}, if one considers $G^{\left(
1\right)  }$ to consist of maps of singletons $\left(  x\rightarrow y\right)
,$ then each of the axioms make quite a bit of sense: the unit is $e\left(
x\right)  =\left(  x\rightarrow x\right)  ,$ the source and range maps are
$s\left(  x\rightarrow y\right)  =x$ and $r\left(  x\rightarrow y\right)  =y,$
associativity ensures composition is okay, and inversion is $\left(
x\rightarrow y\right)  ^{-1}=\left(  y\rightarrow x\right)  .$

\begin{remark}
It might be tempting to replace elements $\gamma$ in $G^{\left(  1\right)  }$
by elements $\left(  s\left(  \gamma\right)  ,r\left(  \gamma\right)  \right)
$ in $G^{\left(  0\right)  }\times G^{\left(  0\right)  }.$ However, often
there will be more than one element of $G^{\left(  1\right)  }$ corresponding
to $\left(  s\left(  \gamma\right)  ,r\left(  \gamma\right)  \right)  .$ See,
e.g., Example \ref{jets}.
\end{remark}

The actual space represented by a groupoid is the orbit space, defined now. We
essentially want the space to be $G^{\left(  0\right)  }$ modulo the
identifications made in $G^{\left(  1\right)  }.$

\begin{definition}
The \emph{orbit} $O_{x}$ of a point $x\in G^{\left(  0\right)  }$ is defined
to be
\[
O_{x}=s\left(  r^{-1}\left(  x\right)  \right)  .
\]
Note that this means that the orbit consists of all points which map to $x$
via an arrow in $G^{\left(  1\right)  }.$ The quotient space $G^{\left(
0\right)  }/\sim,$ where $x\sim y$ if and only if $y\in O_{x}$, is called the
\emph{orbit space}.
\end{definition}

\begin{definition}
A pointed groupoid $\left(  G,O_{x}\right)  $ is a groupoid $G$ together with
a distinguished orbit $O_{x}.$
\end{definition}

Often the orbit space is the actual space we are interested in. In Example
\ref{manifold example} we see that $M$ is the orbit space, and in Example
\ref{quotient example} we see that $X/\Gamma$ is the orbit space.

We want a notion which essentially tells us if the orbit spaces of two
groupoids are the same. For instance, we would like to know that the trivial
groupoid where $G^{\left(  0\right)  }=M$ and $G^{\left(  1\right)  }=M$ is
equivalent to Example \ref{manifold example}. The first guess might be to
define isomorphisms in the categorical sense as appropriate morphisms between
groupoids with inverses. This turns out to be too strong a requirement, so we
introduce a weaker form of equivalence. We begin with the notion of
localization, which is essentially the same procedure that we used to
construct Example \ref{manifold example}. The idea is that we can always take
a groupoid and turn it into a disjoint union of open sets which get identified
via $G^{\left(  1\right)  }.$

\begin{definition}
Let $U=\left\{  U_{i}\right\}  _{i\in I}$ be a cover of $G^{\left(  0\right)
}.$ The \emph{localization} of a groupoid $G$ is the groupoid $G_{U}$ given
by
\[
G_{U}^{\left(  0\right)  }=\coprod\limits_{i\in I}U_{i}=\bigcup\limits_{i\in
I,x\in U_{i}}\left(  i,x\right)
\]
and
\[
G_{U}^{\left(  1\right)  }=\bigcup\limits_{i,j\in I,\gamma\in s^{-1}\left(
U_{i}\right)  \cap r^{-1}\left(  U_{j}\right)  }\left(  i,\gamma,j\right)  .
\]
The unit map is $e\left(  i,x\right)  =\left(  i,e\left(  x\right)  ,i\right)
.$ The source and range maps are $s\left(  i,\gamma,j\right)  =\left(
i,s\left(  \gamma\right)  \right)  ,$ $r\left(  i,\gamma,j\right)  =\left(
j,r\left(  \gamma\right)  \right)  .$ The product is $\left(  i,\gamma
,j\right)  \left(  j,\gamma^{\prime},k\right)  =\left(  i,\gamma\gamma
^{\prime},k\right)  .$
\end{definition}

We have the following definition of equivalence.

\begin{definition}
Two groupoids $G$ and $G^{\prime}$ are equivalent if each has a localization
$G_{U}$ and $G_{U^{\prime}}^{\prime}$ such that $G_{U}$ is isomorphic to
$G_{U^{\prime}}^{\prime}.$
\end{definition}

Note that the property of being equivalent is weaker than the property of
being isomorphic. It will be important to differentiate between equivalence
and isomorphism, since the groupoid structure encodes more than the
equivalence class. In particular, we may consider Ricci flow on trivial
groupoids representing compact manifolds. The limit may not be a manifold, and
hence it is not a trivial groupoid. However, if we consider equivalent
groupoids, there may be a groupoid limit.

\begin{example}
[Localization of trivial groupoid]We see that Example \ref{manifold example}
is a localization of the trivial groupoid, so they are equivalent. Note that
if a Riemannian manifold has a uniform upper bound on sectional curvature
bound, then one can take geodesic balls of a uniform size as the coordinate
patches, as was exploited in \cite{Fuk} and \cite{G1}.
\end{example}

\begin{example}
[Localization of a quotient]If the quotient is a manifold, we see that Example
\ref{quotient example} is equivalent to the trivial groupoid on the quotient
(or orbit space) by taking disjoint copies of the same regularly covered neighborhoods.
\end{example}

Smoothness of a groupoid will allow us to consider the maps in $G^{\left(
1\right)  }$ as smooth diffeomorphisms on some small open sets. In essence,
this makes the maps in $G^{\left(  1\right)  }$ into germs of diffeomorphisms
of $G^{\left(  0\right)  }$. The formal definition is the following.

\begin{definition}
A groupoid $G$ is \emph{smooth} if

\begin{enumerate}
\item $G^{\left(  0\right)  }$ and $G^{\left(  1\right)  }$ are smooth
manifolds (but only assume that $G^{\left(  0\right)  }$ is Hausdorff and
second countable),

\item $e$ is a smooth embedding,

\item $r$ and $s$ are smooth submersions, and

\item multiplication is a smooth map from $\left\{  \left(  \gamma
,\gamma^{\prime}\right)  \in G^{\left(  1\right)  }\times G^{\left(  1\right)
}:s\left(  \gamma\right)  =r\left(  \gamma^{\prime}\right)  \right\}  $ to
$G^{\left(  1\right)  }$ and inversion is a smooth map.
\end{enumerate}
\end{definition}

\begin{remark}
Based on the definition of groupoid, that $r$ is a smooth submersion follows
from $s$ being a smooth submersion, but we include both in the definition to
make it look more symmetric.
\end{remark}

\begin{definition}
We refer to the \emph{dimension} of the groupoid as the dimension of
$G^{\left(  0\right)  }.$
\end{definition}

\begin{definition}
If $r$ and $s$ are local diffeomorphisms, then the groupoid is said to be
\emph{\'{e}tale}. (Note: a map is said to be \'{e}tale if it is a local diffeomorphism.)
\end{definition}

\begin{remark}
Often we will deal with groupoids which are not naturally \'{e}tale. The
groupoid can often be made \'{e}tale by putting the sheaf topology on
$G^{\left(  1\right)  },$ but in general we will not find a need to make our
groupoids \'{e}tale.
\end{remark}

In Examples \ref{manifold example} and \ref{quotient example}, we see that the
arrows come from local diffeomorphisms. This can be made precise with the
following definition.

\begin{definition}
A \emph{local bisection} is a smooth map $\sigma:U\rightarrow G^{\left(
1\right)  },$ where $U\subset G^{\left(  0\right)  }$ is open, such that
$s\circ\sigma$ is the identity and the map $r\circ\sigma:U\rightarrow
r\circ\sigma\left(  U\right)  \ $is a diffeomorphism. We use $\mathcal{B}%
^{loc}\left(  G\right)  $ to refer to the local bisections and $\mathcal{D}%
^{loc}\left(  G\right)  =\left\{  \phi=r\circ\sigma:\sigma\in\mathcal{B}%
^{loc}\left(  G\right)  \right\}  $ to refer to the local diffeomorphisms they generate.
\end{definition}

It is not hard to see that given any element $\gamma\in G^{\left(  1\right)
},$ there is a local bisection $\sigma$ with $\gamma=\sigma\left(  s\left(
\gamma\right)  \right)  $ (see \cite[Prop 5.3]{MM} or \cite[Prop 1.4.9]{Mac}).
The idea is that the property that $r$ and $s$ are submersions is equivalent
to the statement that for any $\gamma\in G^{\left(  1\right)  }$ there is an
open set $U$ containing $\gamma$ such that $\left\{  \left(  s\left(
\gamma^{\prime}\right)  ,r\left(  \gamma^{\prime}\right)  \right)
:\gamma^{\prime}\in U\right\}  $ is the graph of a diffeomorphism. This is
because we may take a local section of $s$ which is transverse to the fibers
of $r,$ giving the graph of a diffeomorphism. Thus we may think of $G^{\left(
1\right)  }$ as containing germs of diffeomorphisms. If the groupoid is
\'{e}tale, then this diffeomorphism is unique up to shrinking the domain and range.

\begin{example}
[Jets of local diffeomorphisms]\label{jets}Given a manifold $M,$ we may define
groupoids of jets of local diffeomorphisms of $M$, denoted $J_{k}=J_{k}\left(
M\right)  ,$ as follows. For each $k$, we define $J_{k}^{\left(  0\right)
}=M.$ We can define $J_{k}^{\left(  1\right)  }$ as pointed diffeomorphisms
$\phi:\left(  U,p\right)  \rightarrow\left(  V,q\right)  $ for open
neighborhoods $U$ of $p$ and $V$ of $q$ modulo an equivalence relation. For
$k=0,$ two maps $\phi:\left(  U,p\right)  \rightarrow\left(  V,q\right)  $ and
$\phi^{\prime}:\left(  U^{\prime},p^{\prime}\right)  \rightarrow\left(
V^{\prime},q^{\prime}\right)  $ are equivalent if $p=p^{\prime}$ and
$q=q^{\prime}.$ For arbitrary $k,$ two maps are equivalent if $p=p^{\prime},$
$q=q^{\prime},$ and all derivatives at $p$ of order less then or equal to $k$
are equal. The source and range maps are defined as $s\left(  \phi\right)  =p$
and $r\left(  \phi\right)  =\phi\left(  p\right)  =q.$

The jet groupoids are not naturally \'{e}tale, though they can be made
\'{e}tale by choosing the sheaf topology. We will generally not do this. Also,
given a Riemannian metric on $M,$ there is a natural Riemannian metric on the
$J_{k}^{\left(  1\right)  }$ defined using the Riemannian metric and the
Riemannian connection. This will be important later, where we use Hausdorff
convergence of closed subsets of $J_{k}^{\left(  1\right)  }$ to define
convergence of Riemannian groupoids.

Given a diffeomorphism $F:M\rightarrow M^{\prime},$ there is a map
$J_{k}^{\left(  1\right)  }\left(  M\right)  $ to $J_{k}^{\left(  1\right)
}\left(  M^{\prime}\right)  $ given by taking $\left[  \gamma\right]  \in
J_{k}^{\left(  1\right)  }\left(  M\right)  $ to $\left[  F\circ\gamma\circ
F^{-1}\right]  \in J_{k}^{\left(  1\right)  }\left(  M^{\prime}\right)  .$ We
denote $\left[  F\circ\gamma\circ F^{-1}\right]  $ by $F_{\ast}\gamma.$

Note that each self-diffeomorphism $F:M\rightarrow M$ induces a global
bisection of each jet groupoid.
\end{example}

A smooth groupoid can be given a Riemannian structure by a putting Riemannian
metric on $G^{\left(  0\right)  }$ that respects the maps in $G^{\left(
1\right)  }$ (i.e., these maps act as isometries).

\begin{definition}
A smooth groupoid $G$ is \emph{Riemannian} if there is a Riemannian metric $g$
on $G^{\left(  0\right)  }$ so that elements of $\mathcal{D}^{loc}\left(
G\right)  $ act as Riemannian isometries.
\end{definition}

There is a natural distance structure on the orbit space of a Riemannian
groupoid (actually, it is only a pseudo-distance, since distinct points may
have zero distance between them).

\begin{definition}
A smooth path $\alpha$ in $G$ is a partition $0=t_{0}\leq t_{1}\leq\cdots\leq
t_{k}=1$ and a sequence
\[
\alpha=\left(  \gamma_{0},\alpha_{1},\gamma_{1},\alpha_{2},\ldots,\alpha
_{k},\gamma_{k}\right)
\]
where $\alpha_{i}:\left[  t_{i-1},t_{i}\right]  \rightarrow G^{\left(
0\right)  }$ is a smooth path and $\gamma_{i}\in G^{\left(  1\right)  }$ with
$\alpha_{i}\left(  t_{i-1}\right)  =r\left(  \gamma_{i-1}\right)  $ and
$\alpha_{i}\left(  t_{i}\right)  =s\left(  \gamma_{i}\right)  .$ The length of
a path is given by
\[
L\left(  \alpha\right)  =\sum_{i=1}^{k}L\left(  \alpha_{i}\right)  .
\]

\end{definition}

\begin{definition}
The pseudometric $d$ on the orbit space of a Riemannian groupoid is given by
\[
d\left(  O_{x},O_{y}\right)  =\inf\left\{  L\left(  \alpha\right)  \right\}
\]
where the infimum is taken over all smooth paths with $s\left(  \gamma
_{0}\right)  =x$ and $r\left(  \gamma_{k}\right)  =y.$ If $d$ is a metric and
the orbits are closed, we say the groupoid is \emph{closed}.
\end{definition}

\begin{remark}
Lott \cite{Lot1} points out that Haefliger \cite{Hae1} and Salem \cite{Sal}
show how to form a closed groupoid by embedding the groupoid into the jet
groupoid $J_{1}$ and taking the closure. In general, the topology on the space
of groupoids will not see the difference between a groupoid and its closure,
much like the Gromov-Hausdorff distance does not see a difference between a
metric space and its completion.
\end{remark}

\begin{definition}
We may define the metric balls $B_{R}\left(  O_{x}\right)  \subset G^{\left(
0\right)  }$ as the union of all orbits which are a distance less than $R$
away from the orbit $O_{x}.$
\end{definition}

Given these definitions, we could define Gromov-Hausdorff distance of
Riemannian groupoids. Instead, we only define $C^{k}$ convergence. The idea is
that we must have local convergence of the Riemannian metrics on $G^{\left(
0\right)  }$ and we must also have local convergence of the arrows $G^{\left(
1\right)  }.$

\begin{definition}
[Convergence of Riemannian groupoids]Let $\left\{  \left(  G_{i}%
,g_{i},O_{x_{i}}\right)  \right\}  _{i=1}^{\infty}$ be a sequence of closed,
pointed, $n$-dimensional Riemannian groupoids and let $\left(  G_{\infty
},g_{\infty},O_{x_{\infty}}\right)  $ be a closed, pointed Riemannian
groupoid. Let $J_{k}$ be the groupoid of $k$-jets of local diffeomorphisms of
$G_{\infty}^{\left(  0\right)  }.$ Then we say that $\lim_{i\rightarrow\infty
}\left(  G_{i},O_{x_{i}}\right)  =\left(  G_{\infty},O_{x_{\infty}}\right)  $
in the pointed $C^{k}$ topology if for all $R>0,$

\begin{enumerate}
\item There exists $I=I\left(  R\right)  $ such that for all $i\geq I$ there
are pointed diffeomorphisms
\[
\phi_{i,R}:B_{R}\left(  O_{x_{\infty}}\right)  \rightarrow B_{R}\left(
O_{x_{i}}\right)
\]
so that
\[
\lim_{i\rightarrow\infty}\phi_{i,R}^{\ast}\left.  g_{i}\right\vert
_{B_{R}\left(  O_{x_{i}}\right)  }=\left.  g_{\infty}\right\vert
_{B_{R}\left(  O_{x_{\infty}}\right)  }%
\]
in $C^{k}\left(  B_{R}\left(  O_{x_{\infty}}\right)  \right)  .$

\item The sets
\[
\phi_{i,R}^{\ast}\left[  s_{i}^{-1}\left(  B_{R/2}\left(  O_{x_{i}}\right)
\right)  \cap r_{i}^{-1}\left(  B_{R}\left(  O_{x_{i}}\right)  \right)
\right]
\]
(see Example \ref{jets}) converge to $s_{\infty}^{-1}\left(  B_{R/2}\left(
O_{x_{i}}\right)  \right)  \cap r_{\infty}^{-1}\left(  B_{R}\left(  O_{x_{i}%
}\right)  \right)  $ in the Hausdorff metric on $J_{k}^{\left(  1\right)
}\left(  G_{\infty}^{\left(  0\right)  }\right)  .$
\end{enumerate}
\end{definition}

\begin{remark}
As noted by Lott \cite{Lot1}, for $k\geq1,$ we need only consider the
convergence in the space of $1$-jets, since the maps are local isometries and
they are entirely determined by their $1$-jets. We keep the $k$ in the
definition here for symmetry in the definition.
\end{remark}

\begin{remark}
In \cite{Fuk} and \cite{G1}, instead of convergence in the space of jets,
convergence in the space of continuous maps is considered. We note that if all
of the arrows can be extended to smooth maps from a fixed domain (such as a
Euclidean ball), then the Arzela-Ascoli theorem tells us that convergence of
the jets implies convergence in $C^{k}$ of the maps. In the examples in the
rest of the paper, the arrows will come from globally defined maps, and we
will therefore deal only with these maps without reference to jets.
\end{remark}

In this paper, we will primarily prove $C^{0}$ convergence. Although it is not
difficult to prove convergence in $C^{\infty},$ we restrict to $C^{0}$ for
clarity of exposition. In all of our examples, the convergence will be
explicit and straightforward.

Groupoid convergence allows one to see collapsing in the following sense.

\begin{definition}
If we start with a sequence of Riemannian groupoids $\left\{  G_{i}\right\}
_{i=1}^{\infty}$ whose orbits are discrete and they converge to a limit
groupoid $G_{\infty}$ such that the orbit space is not discrete, we say that
the sequence \emph{collapses}.
\end{definition}

\section{Solitons on Ricci and Cross Curvature Flows}

Although many of the ideas here apply to higher dimensions (some examples are
given by Lott \cite{Lot1}), we restrict ourselves to dimension 3. In the
sequel, let $\left(  M,g\right)  $ be a three-dimensional Riemannian manifold.

\subsection{Introduction to RF and XCF}

The Ricci flow was first introduced by Hamilton \cite{H1} to study
three-dimensional Riemannian manifolds. The Ricci flow is a solution to the
partial differential equation on Riemannian metrics given by
\begin{equation}
\frac{\partial}{\partial t}g=-2\operatorname{Rc}\left(  g\right)  .
\tag{RF}\label{RF}%
\end{equation}
It is well known that the Ricci flow is weakly parabolic and has a unique
solution for short time (see \cite{H1}). For future use, we note that the
Ricci tensor is invariant under rescaling of the metric, i.e., for any
positive constant $c,$
\begin{equation}
\operatorname{Rc}\left(  cg\right)  =\operatorname{Rc}\left(  g\right)  .
\label{scaling for rc}%
\end{equation}

The cross curvature flow on a three-dimensional manifold was first proposed by
Chow and Hamilton in \cite{CH}. Define the tensor $P^{ij}$ as
\begin{align*}
P^{ij}  &  =R^{ij}-\frac{1}{2}Rg^{ij}\\
&  =g^{ik}g^{j\ell}R_{k\ell}-\frac{1}{2}Rg^{ij}.
\end{align*}
Since we are in dimension $3$, we can diagonalize the Ricci tensor with an
orthonormal frame $\left\{  e_{1},e_{2},e_{3}\right\}  $ and make $P$ diagonal
with $P\left(  \omega^{i},\omega^{i}\right)  $ equal to the sectional
curvature $K\left(  e_{j}\wedge e_{k}\right)  ,$ where $\left\{  \omega
^{1},\omega^{2},\omega^{3}\right\}  $ is the dual coframe and $\left\{
i,j,k\right\}  $ are distinct.

Let $V_{ij}$ be the inverse of $P^{ij}$ (if it exists) and then we define the
cross curvature tensor as
\[
h_{ij}=\left(  \frac{\det P^{ij}}{\det g^{ij}}\right)  V_{ij}.
\]
Notice that $V_{ij}=\frac{1}{\det P^{ij}}\operatorname*{adj}\left(  P\right)
$ so the cross curvature tensor exists even if $P$ is not invertible (though
it may not be an elliptic operator of $g$). We note the following scaling
property of the cross curvature tensor
\begin{equation}
h\left(  cg\right)  =\frac{1}{c}h\left(  g\right)  \label{scaling for h}%
\end{equation}
for a positive constant $c.$

\begin{definition}
The (negative/positive) cross curvature flow ($\pm$XCF) is the flow of
Riemannian metrics solving
\begin{equation}
\frac{\partial}{\partial t}g=-2h\left(  g\right)  \tag{--XCF}%
\end{equation}
or
\begin{equation}
\frac{\partial}{\partial t}g=2h\left(  g\right)  . \tag{+XCF}%
\end{equation}
When we omit the $+$ or $-,$ we are referring to either flow.
\end{definition}

Because our singularity models may change the direction of the flow, it will
often become irrelevant which direction we are considering. However, the
direction is very important for existence results. It was shown that +XCF
exists if the sectional curvature is positive and --XCF exists if the
sectional curvature is negative \cite{Buc}. In other cases, the equation makes
sense, but there may not be a unique solution flow. In the cases of
homogeneous spaces, the partial differential equation reduces to an ordinary
differential equation and thus has a unique solution for a short time.

\subsection{Solitons}

In this section we review soliton techniques for geometric flows. Consider any
geometric flow given by
\begin{equation}
\frac{\partial g}{\partial t}=-2v\left(  g\right)  \label{general equation}%
\end{equation}
where $v$ is a symmetric two tensor which is function of the metric (e.g.,
$\operatorname{Rc},$ $\pm h$). Furthermore, suppose that for any positive
constant $c,$%
\begin{equation}
v\left(  cg\right)  =c^{p}v\left(  g\right)  \label{scaling for v}%
\end{equation}
for some integer $p$ and $v$ is natural, i.e.,
\[
v\left(  \phi^{\ast}g\right)  =\phi^{\ast}\left(  v\left(  g\right)  \right)
\]
for any diffeomorphism $\phi:M\rightarrow M$. If $v=\operatorname{Rc}$ then
$p=0$ by (\ref{scaling for rc}) and if $v=\pm h$ then $p=-1$ by
(\ref{scaling for h}).

\begin{definition}
A \emph{self-similar solution} is a solution of the form
\[
g\left(  t\right)  =\sigma\left(  t\right)  \phi_{t}^{\ast}g_{0}%
\]
where $\sigma$ is a positive function with $\sigma\left(  0\right)  =1$,
$\phi_{t}$ is a one-parameter family of diffeomorphisms of $M$ with $\phi_{0}$
the identity, and $g_{0}$ is a fixed Riemannian metric on $M.$
\end{definition}

\begin{definition}
A \emph{soliton} is a metric $g_{0}$ such that there exists a vector field $X$
on $M$ and a constant $\alpha$ such that
\[
-2v_{0}=L_{X}g_{0}+\alpha g_{0}.
\]
The soliton is said to be a \emph{steady} soliton if $\alpha=0.$
\end{definition}

We give special names to solitons for Ricci flow and cross curvature flow.

\begin{definition}
We refer to solitons of the Ricci flow as \emph{Ricci solitons} and solitons
of the cross curvature flow as \emph{XC solitons}. Note that XC solitons are
solitons for both positive and negative cross curvature flows.
\end{definition}

\begin{remark}
If $p\neq1$ then the metric can be rescaled to produce a soliton with
$\alpha\in\left\{  -1,0,1\right\}  .$
\end{remark}

The following two propositions are the obvious generalizations of \cite[Lemma
2.4 on p. 23]{CK} and \cite[Proposition 1.3 and its successive remarks on pp.
3-4]{RFV2p1}.

\begin{proposition}
[Self-similar iff soliton]\label{self similar iff soliton}If $g\left(
t\right)  $ is a self-similar solution for $t\in\lbrack0,T)$ then $g\left(
0\right)  $ is a soliton. Conversely, if $g_{0}$ is a soliton, then there
exists $T>0$ and a self-similar solution $g\left(  t\right)  $ for
$t\in\lbrack0,T)$ with $g\left(  0\right)  =g_{0}.$
\end{proposition}

\begin{proof}
If $g\left(  t\right)  $ is a self-similar solution, then
\[
g\left(  t\right)  =\sigma\left(  t\right)  \phi_{t}^{\ast}g_{0}.
\]
Differentiating with respect to $t,$ we get%
\[
-2v\left(  \sigma\left(  t\right)  \phi_{t}^{\ast}g_{0}\right)  =\frac
{d\sigma}{dt}\left(  t\right)  \phi_{t}^{\ast}g_{0}+\sigma\left(  t\right)
\phi_{t}^{\ast}\left(  L_{X}g_{0}\right)  ,
\]
where $X$ is the solution to $X\left(  \phi_{t}\left(  p\right)  \right)
=\frac{d}{dt}\phi_{t}\left(  p\right)  $ for all $p\in M.$ By
(\ref{scaling for v}),
\[
v\left(  \sigma\left(  t\right)  \phi_{t}^{\ast}g_{0}\right)  =\left(
\sigma\left(  t\right)  \right)  ^{p}\phi_{t}^{\ast}v\left(  g_{0}\right)
\]
we can drop the pullbacks and get
\begin{equation}
-2\left(  \sigma\left(  t\right)  \right)  ^{p}v\left(  g_{0}\right)
=\frac{d\sigma}{dt}\left(  t\right)  g_{0}+L_{\sigma\left(  t\right)  X}g_{0}.
\label{ladeda}%
\end{equation}
At $t=0,$ this is exactly%
\[
-2v\left(  g_{0}\right)  =\frac{d\sigma}{dt}\left(  0\right)  g_{0}%
+L_{X\left(  0\right)  }g_{0}.
\]

Conversely, if $g\left(  0\right)  $ is a soliton, then
\[
-2v\left(  g_{0}\right)  =L_{X}g_{0}+\alpha g_{0}.
\]
Let $\sigma\left(  t\right)  =\left(  1+\left(  1-p\right)  \alpha t\right)
^{1/\left(  1-p\right)  }$ if $p\neq1$ and $\sigma\left(  t\right)
=\exp\left(  \alpha t\right)  $ if $p=1$ (so $\frac{d\sigma}{dt}=\alpha
\sigma^{p}$). Let $\phi_{t}$ be the diffeomorphisms generated by $\left(
\sigma\left(  t\right)  \right)  ^{p-1}X.$ Then the metric $g\left(  t\right)
=\sigma\left(  t\right)  \phi_{t}^{\ast}g_{0}$ satisfies
\begin{align*}
\frac{\partial}{\partial t}g\left(  t\right)   &  =\frac{d\sigma}{dt}\left(
t\right)  \phi_{t}^{\ast}g_{0}+\phi_{t}^{\ast}\left(  L_{\left(  \sigma\left(
t\right)  \right)  ^{p}X}g_{0}\right) \\
&  =\left(  \sigma\left(  t\right)  \right)  ^{p}\phi_{t}^{\ast}\left(  \alpha
g_{0}+L_{X}g_{0}\right) \\
&  =-2\left(  \sigma\left(  t\right)  \right)  ^{p}\phi_{t}^{\ast}\left(
v\left(  g_{0}\right)  \right) \\
&  =-2v\left(  \sigma\left(  t\right)  \phi_{t}^{\ast}g_{0}\right)  .
\end{align*}

\end{proof}

For the rest of this paper, assume for simplicity that $p\neq1.$ We will
primarily be concerned with $p=0$ for Ricci flow and $p=-1$ for cross
curvature flow. There will always be a corresponding expression if $p=1$ which
we will not provide.

Because of this proposition, we will often interchange the two terms. Note
that we must assume the existence of a solution.

\begin{proposition}
[Canonical form]\label{canonical form prop}Suppose $g\left(  t\right)  $ is a
self similar solution. Then there exist diffeomorphisms $\psi_{t}$ and a
constant $\alpha\in\mathbb{R}$ such that
\[
g\left(  t\right)  =\left(  1+\left(  1-p\right)  \alpha t\right)  ^{1/\left(
1-p\right)  }\psi_{t}^{\ast}g_{0}.
\]

\end{proposition}

\begin{proof}
We have supposed that
\[
g\left(  t\right)  =\sigma\left(  t\right)  \phi_{t}^{\ast}g_{0}.
\]
By (\ref{ladeda}), we have%
\begin{align*}
-2\left(  \sigma\left(  t\right)  \right)  ^{p}v\left(  g_{0}\right)   &
=\frac{d\sigma}{dt}\left(  t\right)  g_{0}+L_{\sigma\left(  t\right)  X}%
g_{0},\\
-2v\left(  g_{0}\right)   &  =\sigma\left(  t\right)  ^{-p}\frac{d\sigma}%
{dt}\left(  t\right)  g_{0}+L_{\sigma\left(  t\right)  ^{1-p}X}g_{0}\\
&  =\frac{1}{1-p}\frac{d\sigma^{1-p}}{dt}g_{0}+L_{\sigma^{1-p}X}g_{0}.
\end{align*}
Differentiating this equation again with respect to $t$ give us
\[
0=\frac{1}{1-p}\frac{d^{2}\sigma^{1-p}}{d^{2}t}g_{0}+L_{\tilde{X}}g_{0},
\]
where $\tilde{X}\left(  t\right)  =\frac{d\sigma^{1-p}}{dt}X+\sigma^{1-p}%
\frac{dX}{dt}.$ So either $\frac{d^{2}\sigma^{1-p}}{dt^{2}}=0$ or $g_{0}%
=L_{Y}g_{0}$ with $Y=-\tilde{X}/\left(  \frac{1}{1-p}\frac{d^{2}\sigma^{1-p}%
}{d^{2}t}\right)  .$ In the first case,
\[
\sigma\left(  t\right)  =\left(  1+\left(  1-p\right)  \alpha t\right)
^{1/\left(  1-p\right)  },
\]
since $g\left(  0\right)  =g_{0}.$ In the second case,
\[
-2v\left(  g_{0}\right)  =L_{\beta X+\gamma Y}g_{0},
\]
where $\beta=\sigma^{1-p}$ and $\gamma=\left(  \frac{1}{1-p}\frac
{d\sigma^{1-p}}{dt}\right)  $ and so we may choose $\alpha=0.$ The proof is
completed by Proposition \ref{self similar iff soliton}.
\end{proof}

\begin{corollary}
All Ricci solitons can be put in the form
\[
g\left(  t\right)  =\left(  1+\alpha t\right)  \psi_{t}^{\ast}g_{0},
\]
and all XC solitons can be put in the form
\[
g\left(  t\right)  =\left(  1+2\alpha t\right)  ^{1/2}\psi_{t}^{\ast}g_{0}.
\]

\end{corollary}

\begin{remark}
Following Lott \cite{Lot1}, we will base our self-similar solutions at
$g_{1}=g\left(  1\right)  $ instead of $g_{0}=g\left(  0\right)  .$ In this
case, the canonical forms are $g\left(  t\right)  =\alpha t\psi_{t}^{\ast
}g_{1}$ for Ricci flow and $g\left(  t\right)  =\left(  2\alpha t\right)
^{1/2}\psi_{t}^{\ast}g_{1}$ for cross curvature flow.
\end{remark}

\subsection{Theory of singularities}

In order to understand the geometry of a limit solution, one must look at the
appropriate length scale. For instance, given any Riemannian manifold $\left(
\mathcal{M},g\right)  $ and a point $p\in\mathcal{M},$ one could consider the
manifold gotten by the limit of $\left(  \mathcal{M}^{n},sg\right)  $ where
$s\rightarrow\infty.$ Since the space is a Riemannian manifold, this will
converge in the Gromov-Hausdorff sense to Euclidean space $\left(
\mathbb{R}^{n},g_{\mathbb{E}}\right)  .$ On the other hand, if one takes the
unit sphere metric $\left(  S^{n},g_{S^{n}\left(  1\right)  }\right)  $ and
looks at the limit $\left(  S^{n},sg_{S^{n}\left(  1\right)  }\right)  $ where
$s\rightarrow0,$ it is clear that the sectional curvatures go to infinity. As
we are taking limits, in order to understand the geometry of a particular
solution, we will wish to rescale in such a way that we get a reasonable limit
that has, if possible, nonzero curvatures. This is what we will call the
\emph{geometric limit}. The general process for rescaling is to rescale so
that the maximum sectional curvature (in absolute value) does not go to zero
or infinity. If $g\left(  t\right)  $ is a solution to the geometric equation
(\ref{general equation}), we will do a parabolic rescaling so that the limit
is also a solution to the flow. Suppose $g\left(  t\right)  $ is defined on a
maximal time interval $[0,T).$ The usual rescaling as the flow goes to the
singularity at $T\in(0,\infty]$ is as follows (see \cite{H2}, \cite{CK}
\cite{CLN}). We take an increasing sequence $t_{i}$ converging to $T$ (or
going to infinity if $T=\infty$) and consider the sequence of metrics%
\[
g_{i}\left(  t\right)  =M\left(  t_{i}\right)  g\left(  \frac{t}{M\left(
t_{i}\right)  }+t_{i}\right)  .
\]
These solutions have the property that $g_{i}\left(  0\right)  =g\left(
t_{i}\right)  $ for some function $M\left(  t\right)  .$ Note that if we take
\begin{equation}
M\left(  t\right)  =\sup_{\mathcal{M}}\left\vert \operatorname{Rm}\left(
g\left(  t\right)  \right)  \right\vert \label{M def}%
\end{equation}
then $\sup_{\mathcal{M}}\left\vert \operatorname{Rm}\left(  g_{i}\left(
0\right)  \right)  \right\vert =1.$

In following Lott \cite{Lot1}, we consider a continuous deformation (so
instead of taking a sequence $t_{i},$ we take a parameter $s$) with base
metric $g_{s}\left(  1\right)  $ (instead of $g_{i}\left(  0\right)  $). If
$T=\infty,$ we will consider deformations of the form
\begin{equation}
g_{s}\left(  t\right)  =f\left(  s\right)  g\left(  \frac{t-1}{\left(
f\left(  s\right)  \right)  ^{1-p}}+s\right)  , \label{g_s definition}%
\end{equation}
where $p$ is defined by (\ref{scaling for v}). These have the property that
$g_{s}\left(  1\right)  =g\left(  s\right)  $ and that $g_{s}\left(  t\right)
$ is defined for $t$ in%
\[
\lbrack1-\left(  f\left(  s\right)  \right)  ^{1-p}s,\infty).
\]
We also see that $g_{s}\left(  t\right)  $ satisfies
\begin{align*}
\frac{\partial}{\partial t}g_{s}  &  =-2\left(  f\left(  s\right)  \right)
^{p}v\left(  g\left(  \frac{t-1}{\left(  f\left(  s\right)  \right)  ^{1-p}%
}+s\right)  \right) \\
&  =-2v\left(  g_{s}\right)  .
\end{align*}
We will then consider limits as $s\rightarrow\infty,$
\[
g_{\infty}\left(  t\right)  =\lim_{s\rightarrow\infty}\phi_{s}^{\ast}%
g_{s}\left(  t\right)
\]
where $\phi_{s}$ are appropriately chosen diffeomorphisms. Note that if
$f\left(  s\right)  =s^{-1/\left(  1-p\right)  }$ then $g_{s}\left(  t\right)
=s^{-1/\left(  1-p\right)  }g\left(  st\right)  ,$ and
\[
g_{\infty}\left(  t\right)  =\lim_{s\rightarrow\infty}s^{-1/\left(
1-p\right)  }\phi_{s}^{\ast}g\left(  st\right)  ,
\]
which will be a common rescaling (note that for Ricci flow we have $-1/\left(
1-p\right)  =-1$ and for cross curvature flow we have $-1/\left(  1-p\right)
=-1/2$). In this case, the limit will be defined for $t$ in $[0,\infty).$

In the case of $T<\infty,$ we will instead look at limits defined by
\begin{equation}
g_{s}\left(  t\right)  =f\left(  s\right)  g\left(  T-\left(  \frac
{t-1}{\left(  f\left(  s\right)  \right)  ^{1-p}}+s\right)  \right)  .
\label{other g_s definition}%
\end{equation}
These have the property that $g_{s}\left(  1\right)  =g\left(  T-s\right)  ,$
$g_{s}\left(  t\right)  $ is defined for $t$ in the interval
\[
(1-\left(  f\left(  s\right)  \right)  ^{1-p}s,~\left(  f\left(  s\right)
\right)  ^{1-p}\left(  T-1-s\right)  +1],
\]
and
\begin{equation}
\frac{\partial}{\partial t}g_{s}=2v\left(  g_{s}\right)  \label{backward flow}%
\end{equation}
(notice that the sign is flipped). In this case, to look at the solution near
the singularity, we look at the limit
\[
g_{T}\left(  t\right)  =\lim_{s\rightarrow0}\phi_{s}^{\ast}g_{s}\left(
t\right)
\]
for some choice of diffeomorphisms $\phi_{s}.$ Note that if $f\left(
s\right)  =s^{-1/\left(  1-p\right)  }$ then the limit is defined for $t$ in
$\left(  0,\infty\right)  $ and the limit looks like%
\[
g_{T}\left(  t\right)  =\lim_{s\rightarrow0}s^{-1/\left(  1-p\right)  }%
\phi_{s}^{\ast}g_{s}\left(  st\right)  .
\]

It is often more important for us to understand how a particular solution of a
flow compares with other solutions. In this case, we will consider certain
classes of singularities. As introduced by Hamilton \cite{H2}, one can
separate solutions into 4 classes of solutions (we take $M\left(  t\right)  $
defined by (\ref{M def})):

\begin{itemize}
\item[Type I.] $T<\infty$ and $\sup\left(  T-t\right)  ^{1/\left(  1-p\right)
}M\left(  t\right)  <\infty$

\item[Type IIa.] $T<\infty$ and $\sup\left(  T-t\right)  ^{1/\left(
1-p\right)  }M\left(  t\right)  =\infty$

\item[Type IIb.] $T=\infty$ and $\sup t^{1/\left(  1-p\right)  }M\left(
t\right)  =\infty$

\item[Type III.] $T=\infty$ and $\sup t^{1/\left(  1-p\right)  }M\left(
t\right)  <\infty.$
\end{itemize}

This singularity theory gives a canonical rescaling factor of $t^{1/\left(
1-p\right)  }$ designed to give limit soliton metrics based on the canonical
form of soliton metrics described in Proposition \ref{canonical form prop}. It
is significant that this rescaling is chosen by the flow and not by the
solution itself.

\begin{remark}
Another canonical rescaling one might propose is one such that the volume is
unchanged (often called normalized Ricci flow and normalized cross curvature
flow), as used in \cite{H1}, \cite{IJ}, \cite{CNS}, and others. We do not
treat this particular rescaling, arguing that the geometric rescaling that
keeps the curvatures bounded and the singularity rescaling are more natural in
most of the cases we give here. In many of the cases we treat, rescaling so
that volume is unchanged will not prevent collapsing and convergence will
usually be to a collapsed flat manifold.
\end{remark}

For Ricci flow on three-dimensional manifolds, it is an interesting fact that
most of the homogeneous solutions are Type III, with the exception of $S^{3}$
and $S^{2}\times\mathbb{R}$. For negative cross curvature flow, we will
actually find a Type IIb solution.

When looking at a Type III solution, we will look at the \emph{Type III limit
solution} described by the limit of \emph{Type III rescalings}
\[
s^{-1/\left(  1-p\right)  }g\left(  st\right)
\]
as $s\rightarrow\infty.$ When looking at a Type I solution, we will look at
the \emph{Type I limit solution} described by the limit of \emph{Type I
rescalings}
\[
s^{-1/\left(  1-p\right)  }g\left(  T-st\right)
\]
as $s\rightarrow0.$ Notice that the negative sign makes this a solution not of
the original flow, but of the backward flow (\ref{backward flow}). For Type II
solutions, we will need to find an appropriate geometric rescaling, since the
rescalings we have proposed will take the sectional curvature to infinity.

The existence of limit solutions described in the previous paragraph is not
guaranteed. By Hamilton's compactness theorem \cite{H3}, if we had a uniform
lower bound on the injectivity radius of the scaled solutions, one could take
this limit for Ricci flow (and it is not hard to construct a similar theorem
for XCF). In the Type I Ricci flow case, Perelman showed that this bound
exists \cite{P1}. For Type III Ricci flow this is not a reasonable assumption
since collapsing does happen. Similarly, collapsing can occur for XCF.
However, the idea of limit solutions can apply in the setting of Riemannian
groupoids, and this is where the limits will be taken. Instead of Hamilton's
theorem, Lott's compactness theorem for Ricci flow on Riemannian groupoids
(\cite{Lot1}, see also \cite{Fuk} and \cite{G1} for some of the geometric
ideas) may be used to extract a limit.

\section{3D homogeneous solutions}

In this section we review the results on Ricci flow on three-dimensional
homogeneous geometries and give the results on cross curvature flow. Solutions
of the Ricci flow on three-dimensional, simply connected, homogeneous
geometries were first described by Isenberg-Jackson \cite{IJ} (see also
\cite{KM}, \cite[Chapter 1]{CK}, \cite[Chapter 4, Section 7]{CLN}). In
general, we may start with a basis of left-invariant vector fields
$F_{1},F_{2},F_{3}$ and consider the class of left invariant metrics such that
this frame is orthogonal (but the length of the vectors in the frame is
arbitrary). Solutions of the negative cross curvature flow on simply connected
homogeneous geometries were first described by Cao-Ni-Saloff-Coste \cite{CNS}.
The homogeneous expanding solitons on $\operatorname{Nil}$ and
$\operatorname{Sol}$ were described by Baird-Danielo \cite{BD} and Lott
\cite{Lot1}. The results below on Ricci flow are due to Lott \cite{Lot1}; in
some cases we give different coordinate representations in an attempt to make
the limit groupoids especially clear. The results on cross curvature flow are new.

\begin{remark}
We choose to follow the convention of Lott \cite{Lot1} that the brackets of
the frame look like $\left[  F_{i},F_{j}\right]  =c_{ij}^{k}F_{j}$ where
$c_{ij}^{k}$ are in $\left\{  -1,0,1\right\}  .$ This is different from the
conventions in \cite{KM} and \cite{CNS}, and so curvatures may look slightly
different because of the discrepancy in the definition of $A,B,C.$ We will try
to point out the discrepancies in each example.
\end{remark}

The goal of this section is to find the limits of collapsing solutions of
Ricci flow and cross curvature flow. In the process, we find Ricci and cross
curvature solitons which occur in the limit. The process is as follows. First,
look at the asymptotic solutions of the flow on simply connected geometries.
Since each of the following simply connected geometries is diffeomorphic to
$\mathbb{R}^{3},$ there is a wealth of diffeomorphisms available, primarily
rescaling of the coordinates. We need to find diffeomorphisms so that the
metrics in coordinates which are pulled back by these diffeomorphisms do not
degenerate. The limit geometry may be the same, in which case the geometry
admits a soliton metric, or it may be different, in which case the geometry
converges to another geometry. Lastly, for a compact homogenous manifold, we
consider the equivalent groupoid consisting of the universal cover together
with arrows described by the action of the fundamental group as deck
transformations. If the limit of the arrows describes a continuous group
$\mathcal{D}^{loc}$, there is collapsing. We note that it is extremely
important to find the right coordinates so that the diffeomorphisms may be
written out explicitly.

We shall only look at four of the possible three-dimensional homogeneous
geometries because we wish to emphasize the importance of the changing
diffeomorphisms. One may consider the remaining geometries, but the effects of
the diffeomorphisms is much more trivial than the effects in
$\operatorname{Nil},$ $\operatorname{Sol},$ $\widetilde{\operatorname{SL}%
}\left(  2,\mathbb{R}\right)  ,$ and $\widetilde{\operatorname{Isom}}\left(
\mathbb{E}^{2}\right)  .$

\subsection{$\operatorname{Nil}$}

Recall that $\operatorname{Nil}$ consists of the unit upper triangular
matrices,
\[
\left(
\begin{array}
[c]{ccc}%
1 & a & c\\
0 & 1 & b\\
0 & 0 & 1
\end{array}
\right)  \left(
\begin{array}
[c]{ccc}%
1 & x & z\\
0 & 1 & y\\
0 & 0 & 1
\end{array}
\right)  ,
\]
and so the group multiplication is
\[
\left(  a,b,c\right)  \left(  x,y,z\right)  =\left(  x+a,y+b,z+c+ay\right)  .
\]
We easily see that the following global vector fields are left invariant%
\[
F_{1}=\frac{\partial}{\partial z},\qquad F_{2}=\frac{\partial}{\partial
x},\qquad F_{3}=\frac{\partial}{\partial y}+x\frac{\partial}{\partial z},
\]
and we easily see that
\[
\left[  F_{2},F_{3}\right]  =F_{1}%
\]
and all other brackets are zero. We can also define the dual forms as
\[
\theta_{1}=dz-xdy,\qquad\theta_{2}=dy,\qquad\theta_{3}=dx.
\]
It is then clear that the following metrics are all left invariant%
\begin{align}
g  &  =A\theta_{1}^{2}+B\theta_{2}^{2}+C\theta_{3}^{2}\nonumber\\
&  =A\left(  dz-xdy\right)  ^{2}+Bdy^{2}+Cdx^{2}. \label{Nil metric}%
\end{align}
Note that by appropriate scaling of the coordinates, we see that this is
actually only a one-parameter family of metrics up to diffeomorphism. It is
not difficult to see that the sectional curvatures for these metrics are%
\begin{align*}
K\left(  F_{2}\wedge F_{3}\right)   &  =-\frac{3A}{4BC},\\
K\left(  F_{3}\wedge F_{1}\right)   &  =\frac{A}{4BC},\\
K\left(  F_{1}\wedge F_{2}\right)   &  =\frac{A}{4BC}.
\end{align*}
Note that this $F_{1}$ is half that used in \cite{KM} and \cite{CNS}, so our
$A$ is 1/4 the corresponding coefficient in those papers.

\subsubsection{Ricci Flow}

It is well known that the Ricci flow on the metric (\ref{Nil metric}) has the
form
\begin{align*}
\frac{dA}{dt}  &  =-\frac{A^{2}}{BC},\\
\frac{dB}{dt}  &  =\frac{A}{C},\\
\frac{dC}{dt}  &  =\frac{A}{B},
\end{align*}
(see \cite{IJ} \cite{KM}). The solution is
\begin{align*}
A\left(  t\right)   &  =A_{0}\left(  3\frac{A_{0}}{B_{0}C_{0}}t+1\right)
^{-1/3},\\
B\left(  t\right)   &  =B_{0}\left(  3\frac{A_{0}}{B_{0}C_{0}}t+1\right)
^{1/3},\\
C\left(  t\right)   &  =C_{0}\left(  3\frac{A_{0}}{B_{0}C_{0}}t+1\right)
^{1/3}.
\end{align*}

Notice that the sectional curvatures all behave like $t^{-1},$ so the solution
is Type III. We may pull the metric back by the diffeomorphisms%
\[
\phi_{t}\left(  x,y,z\right)  =\left(  \frac{x}{t^{1/6}},\frac{y}{t^{1/6}%
},t^{1/6}z\right)
\]
to get the metrics
\begin{align*}
\phi_{t}^{\ast}g\left(  t\right)   &  =A_{0}\left(  3\frac{A_{0}}{B_{0}C_{0}%
}+\frac{1}{t}\right)  ^{-1/3}\left(  dz-t^{-1/2}xdy\right)  ^{2}+B_{0}\left(
3\frac{A_{0}}{B_{0}C_{0}}+\frac{1}{t}\right)  ^{1/3}dy^{2}\\
&  \;\;\;\;\;+C_{0}\left(  3\frac{A_{0}}{B_{0}C_{0}}+\frac{1}{t}\right)
^{1/3}dx^{2}.
\end{align*}
We note that as $t\rightarrow\infty$, this converges to the Euclidean metric.
However, this is because the Ricci flow causes the metric to spread out. We
may counteract this by rescaling the metric via a Type III rescaling. So,
instead, we consider
\[
g_{s}\left(  t\right)  =\frac{1}{s}g\left(  st\right)
\]
where $s\rightarrow\infty.$ The idea is that $g_{s}\left(  1\right)  $ is the
long term behavior of $g\left(  t\right)  $ after rescaling. We pull back by
different diffeomorphisms%
\[
\psi_{s}\left(  x,y,z\right)  =\left(  s^{1/3}x,s^{1/3}y,s^{2/3}z\right)
\]
and get%
\begin{align*}
\frac{1}{s}\psi_{s}^{\ast}g_{s}  &  =A_{0}\left(  3\frac{A_{0}}{B_{0}C_{0}%
}t+\frac{1}{s}\right)  ^{-1/3}\left(  dz-xdy\right)  ^{2}+B_{0}\left(
3\frac{A_{0}}{B_{0}C_{0}}t+\frac{1}{s}\right)  ^{1/3}dy^{2}\\
&  \;\;\;\;\;+C_{0}\left(  3\frac{A_{0}}{B_{0}C_{0}}t+\frac{1}{s}\right)
^{1/3}dx^{2}.
\end{align*}
As $s\rightarrow\infty,$ we get the limit Ricci flow
\begin{align*}
g_{\infty}\left(  t\right)   &  =A_{0}\left(  3\frac{A_{0}}{B_{0}C_{0}%
}t\right)  ^{-1/3}\left(  dz-xdy\right)  ^{2}+B_{0}\left(  3\frac{A_{0}}%
{B_{0}C_{0}}t\right)  ^{1/3}dy^{2}\\
&  \;\;\;\;\;+C_{0}\left(  3\frac{A_{0}}{B_{0}C_{0}}t\right)  ^{1/3}dx^{2}.
\end{align*}
Notice that if we pull back by the diffeomorphism
\[
\tilde{\psi}\left(  x,y,z\right)  =\left(  C_{0}^{-1/2}\left(  \frac{3A_{0}%
}{B_{0}C_{0}}\right)  ^{-1/6}x,B_{0}^{-1/2}\left(  \frac{3A_{0}}{B_{0}C_{0}%
}\right)  ^{-1/6}y,B_{0}^{-1/2}C_{0}^{-1/2}\left(  \frac{3A_{0}}{B_{0}C_{0}%
}\right)  ^{-1/3}z\right)
\]
we get
\begin{align*}
\tilde{g}_{\infty}\left(  t\right)   &  =\tilde{\psi}^{\ast}g_{\infty}\left(
t\right) \\
&  =\frac{1}{t^{1/3}}\left(  dz-xdy\right)  ^{2}+t^{1/3}dy^{2}+t^{1/3}dx^{2}.
\end{align*}
This is the $\operatorname{Nil}$ soliton from \cite{Lot1} and \cite{BD}. We
easily see that
\[
\tilde{g}_{\infty}\left(  t\right)  =t\phi_{t}^{\ast}g
\]
where
\[
g=\left(  dz-xdy\right)  ^{2}+dy^{2}+dx^{2}%
\]
and
\[
\phi_{t}\left(  x,y,z\right)  =\left(  t^{-1/3}x,t^{-1/3}y,t^{-2/3}z\right)
.
\]

Now consider compact quotients of $\operatorname{Nil}.$ We can also compute
the limit of isometries on $g\left(  t\right)  $, which are
\[
\psi_{s}^{-1}\circ\gamma_{k,\ell,m}\circ\psi_{s}\left(  x,y,z\right)  =\left(
x+s^{-1/3}k,y+s^{-1/3}\ell,z+s^{-2/3}m+s^{-1/3}ky\right)  .
\]
If we take the limit $s\rightarrow\infty,$ the the group converges to the
following isometries of $g_{\infty}\left(  t\right)  $%
\[
\gamma_{u,v,w}\left(  x,y,z\right)  =\left(  x+u,y+v,z+w+uy\right)  ,
\]
where $u,v,w$ are real numbers gotten by choosing numbers $k\left(  s\right)
,$ $\ell\left(  s\right)  ,$ $m\left(  s\right)  $ such that
\begin{align*}
\lim_{s\rightarrow\infty}s^{-1/3}k\left(  s\right)   &  =u,\\
\lim_{s\rightarrow\infty}s^{-1/3}\ell\left(  s\right)   &  =v,\\
\lim_{s\rightarrow\infty}s^{-2/3}m\left(  s\right)   &  =w.
\end{align*}
Note that this means that $k,\ell,m$ must become large. Now, suppose we start
with a compact quotient of $\operatorname{Nil}.$ Then the fundamental group
may be represented as a subgroup of $\operatorname{Nil}$ acting freely,
properly discontinuously. Still, it is possible for $k,\ell,m$ to become
infinitely large, since the lattice must extend through the entirety of
$\operatorname{Nil}$ for the quotient to be compact. Thus the lattice
converges to all of the group elements, and thus the renormalized manifold
converges to a groupoid whose orbit space is a point. It is interesting to
note that this is, in some sense, the optimal geometry for $\operatorname{Nil}%
,$ as Gromov's almost flat manifold theorem states that any almost flat
manifold must be an infranilmanifold \cite{Gr} (see also \cite{BK}).

\subsubsection{Cross Curvature Flow}

From \cite{CNS}, we see that the negative cross curvature flow on
$\operatorname{Nil}$ has solutions%
\begin{align*}
A\left(  t\right)   &  =A_{0}\left(  1+7R_{0}^{2}t\right)  ^{-1/14},\\
B\left(  t\right)   &  =B_{0}\left(  1+7R_{0}^{2}t\right)  ^{3/14},\\
C\left(  t\right)   &  =C_{0}\left(  1+7R_{0}^{2}t\right)  ^{3/14},
\end{align*}
where $R_{0}=-A_{0}/\left(  2B_{0}C_{0}\right)  $ (this is one fourth the
value in \cite{CNS}). Notice that all sectional curvatures behave
asymptotically like $t^{-1/2},$ so there may be a XC soliton metric on
$\operatorname{Nil}$. We can consider%
\begin{align*}
\frac{1}{s^{1/2}}\phi_{s}^{\ast}g\left(  st\right)   &  =A_{0}\frac{\left(
1+7R_{0}^{2}st\right)  ^{-1/14}s^{4/7}}{s^{1/2}}\left(  dz-xdy\right)
^{2}+B_{0}\frac{\left(  1+7R_{0}^{2}st\right)  ^{3/14}s^{2/7}}{s^{1/2}}%
dy^{2}\\
&  \;\;\;\;\;+C_{0}\frac{\left(  1+7R_{0}^{2}st\right)  ^{3/14}s^{2/7}%
}{s^{1/2}}dx^{2}%
\end{align*}
whose limit as $s\rightarrow\infty$ is
\[
A_{0}\left(  7R_{0}^{2}t\right)  ^{-1/14}\left(  dz-xdy\right)  ^{2}%
+B_{0}\left(  7R_{0}^{2}t\right)  ^{3/14}dy^{2}+C_{0}\left(  7R_{0}%
^{2}t\right)  ^{3/14}dx^{2},
\]
where
\[
\phi_{s}\left(  x,y,z\right)  =\left(  s^{1/7}x,s^{1/7}y,s^{2/7}z\right)  .
\]
Note that we may pull back by
\[
\tilde{\psi}\left(  x,y,z\right)  =\left(  C_{0}^{-1/2}\left(  7R_{0}%
^{2}\right)  ^{-3/28}x,B_{0}^{-1/2}\left(  7R_{0}^{2}\right)  ^{-3/28}%
y,\left(  B_{0}C_{0}\right)  ^{-1/2}\left(  7R_{0}^{2}\right)  ^{-3/14}%
z\right)
\]
to get%
\[
g\left(  t\right)  =\frac{2}{\sqrt{7}t^{1/14}}\left(  dz-xdy\right)
^{2}+t^{3/14}\left(  dy^{2}+dx^{2}\right)  .
\]
This is a XC soliton, since%
\[
g\left(  t\right)  =t^{1/2}\psi_{t}^{\ast}\left(  \frac{2}{\sqrt{7}}\left(
dz-xdy\right)  ^{2}+dy^{2}+dx^{2}\right)
\]
where
\[
\psi_{t}\left(  x,y,z\right)  =\left(  t^{-1/7}x,t^{-1/7}y,t^{-2/7}z\right)
.
\]

We can also compute the limit of isometries on $g\left(  t\right)  $, which
are
\[
\phi_{s}^{-1}\circ\gamma_{k,\ell,m}\circ\phi_{s}\left(  x,y,z\right)  =\left(
x+s^{-1/7}k,y+s^{-1/7}\ell,z+s^{-2/7}m+s^{-1/7}ky\right)  .
\]
The group limit looks much like that for the Ricci flow, in that we get group
elements%
\[
\gamma_{u,v,w}=\left(  x+u,y+v,z+w+uy\right)
\]
where $u,v,w$ are real numbers gotten by choosing numbers $k\left(  s\right)
,$ $\ell\left(  s\right)  ,$ $m\left(  s\right)  $ such that
\begin{align*}
\lim_{s\rightarrow\infty}s^{-1/7}k\left(  s\right)   &  =u,\\
\lim_{s\rightarrow\infty}s^{-1/7}\ell\left(  s\right)   &  =v,\\
\lim_{s\rightarrow\infty}s^{-2/7}m\left(  s\right)   &  =w.
\end{align*}
Once again, we get convergence to a Riemannian groupoid whose orbit space is a point.

\subsection{$\operatorname{Sol}$}

The group $\operatorname{Sol}$ is a Lie group on $\mathbb{R}^{3}$ with a group
action given by
\[
\left(  a,b,c\right)  \left(  x,y,z\right)  =\left(  e^{-c}x+a,e^{c}%
y+b,z+c\right)  .
\]
One can easily see that the frame
\[
F_{1}=e^{-z}\frac{\partial}{\partial x}+e^{z}\frac{\partial}{\partial
y},\;\;\;F_{2}=-\frac{\partial}{\partial z},\;\;\;F_{3}=e^{-z}\frac{\partial
}{\partial x}-e^{z}\frac{\partial}{\partial y}%
\]
is left invariant and satisfies
\begin{align*}
\left[  F_{1},F_{2}\right]   &  =-F_{3},\\
\left[  F_{2},F_{3}\right]   &  =F_{1},\\
\left[  F_{3},F_{1}\right]   &  =0.
\end{align*}
We have a family of left invariant metrics given by%
\[
g=A\left(  e^{z}dx+e^{-z}dy\right)  ^{2}+Bdz^{2}+C\left(  e^{z}dx-e^{-z}%
dy\right)  ^{2}.
\]
This is really only a two-parameter family up to diffeomorphism, since we may
rescale $x$ and $y$, making $A$ and $C$ only well-defined up to their ratio.
We may also find it useful to use alternate coordinates, which give%

\begin{equation}
g=\left(  d\tilde{z}-2\sqrt{\frac{A}{BC}}\tilde{x}d\tilde{y}\right)
^{2}+d\tilde{y}^{2}+\left(  d\tilde{x}-2\sqrt{\frac{C}{AB}}\tilde{z}d\tilde
{y}\right)  ^{2}, \label{alt sol metric}%
\end{equation}
by the map described by%
\begin{align*}
x  &  =e^{\frac{-\tilde{y}}{\sqrt{B}}}\left(  \frac{\tilde{x}}{\sqrt{C}}%
+\frac{\tilde{z}}{\sqrt{A}}\right)  ,\\
y  &  =e^{\frac{\tilde{y}}{\sqrt{B}}}\left(  -\frac{\tilde{x}}{\sqrt{C}}%
+\frac{\tilde{z}}{\sqrt{A}}\right)  ,\\
z  &  =\frac{\tilde{y}}{\sqrt{B}}%
\end{align*}
with inverse%
\begin{align*}
\tilde{x}  &  =\frac{1}{2}\sqrt{C}\left(  e^{z}x-e^{-z}y\right)  ,\\
\tilde{y}  &  =\sqrt{B}z,\\
\tilde{z}  &  =\frac{1}{2}\sqrt{A}\left(  e^{z}x+e^{-z}y\right)  .
\end{align*}
It is not hard to see that the sectional curvatures are%
\begin{align*}
K\left(  F_{2}\wedge F_{3}\right)   &  =\frac{\left(  A-C\right)  ^{2}-4A^{2}%
}{4ABC},\\
K\left(  F_{3}\wedge F_{1}\right)   &  =\frac{\left(  A+C\right)  ^{2}}%
{4ABC},\\
K\left(  F_{1}\wedge F_{2}\right)   &  =\frac{\left(  A-C\right)  ^{2}-4C^{2}%
}{4ABC}.
\end{align*}

An isometry $\gamma_{\left(  a,b,c\right)  }\left(  x,y,z\right)  =\left(
a,b,c\right)  \left(  x,y,z\right)  $ (expressed in the first coordinate
chart) can be brought back to an action on $\left(  \tilde{x},\tilde{y}%
,\tilde{z}\right)  $ as%
\begin{align}
&  \gamma_{\left(  a,b,c\right)  }\left(  \tilde{x},\tilde{y},\tilde{z}\right)
\label{sol isometries}\\
&  =\left(  \tilde{x}+e^{\tilde{y}/\sqrt{B}+c}a-e^{-\tilde{y}/\sqrt{B}%
-c}b,~\tilde{y}+\sqrt{B}c,~\tilde{z}+e^{\tilde{y}/\sqrt{B}+c}a+e^{-\tilde
{y}/\sqrt{B}-c}b\right)  .\nonumber
\end{align}

Note that our $F_{2}$ is one half that in \cite{KM}, and so our $B$ is 1/4 the
corresponding $B$ in that paper and \cite{CNS}, but agrees with \cite{Lot1}.

\subsubsection{Ricci Flow}

From \cite{KM} we see that
\begin{align*}
A,C  &  \sim\sqrt{A_{0}C_{0}},\\
B  &  \sim4t,\\
A-C  &  \sim\frac{E}{t}.
\end{align*}
The sectional curvatures look like
\begin{align*}
K\left(  F_{2}\wedge F_{3}\right)   &  \sim\frac{E_{1}}{t},\\
K\left(  F_{3}\wedge F_{1}\right)   &  \sim\frac{E_{2}}{t},\\
K\left(  F_{1}\wedge F_{2}\right)   &  \sim\frac{E_{3}}{t},
\end{align*}
for some constants $E_{1},E_{2},E_{3}.$ Thus the solution is Type III. We may
take the Type III limit as
\begin{align*}
&  \lim_{s\rightarrow\infty}s^{-1}\phi_{s}^{\ast}g\left(  st\right) \\
&  =\lim_{s\rightarrow\infty}\left[  s^{-1}\sqrt{A_{0}C_{0}}\left(
e^{z}s^{1/2}dx+e^{-z}s^{1/2}dy\right)  ^{2}+4tdz^{2}+s^{-1}\sqrt{A_{0}C_{0}%
}\left(  e^{z}s^{1/2}dx-e^{-z}s^{1/2}dy\right)  ^{2}\right] \\
&  =\sqrt{A_{0}C_{0}}\left(  e^{z}dx+e^{-z}dy\right)  ^{2}+4tdz^{2}%
+\sqrt{A_{0}C_{0}}\left(  e^{z}dx-e^{-z}dy\right)  ^{2}\\
&  =\sqrt{A_{0}C_{0}}\left(  e^{2z}dx^{2}+e^{-2z}dy^{2}\right)  +4tdz^{2},
\end{align*}
where $\phi_{s}\left(  x,y,z\right)  =\left(  s^{1/2}x,s^{1/2}y,z\right)  .$
We can pull back by
\[
\psi\left(  x,y,z\right)  =\left(  \left(  A_{0}C_{0}\right)  ^{-1/4}x,\left(
A_{0}C_{0}\right)  ^{-1/4}y,z\right)
\]
to get the limit soliton
\[
e^{2z}dx^{2}+e^{-2z}dy^{2}+4tdz^{2}.
\]
This is the $\operatorname{Sol}$ soliton described in \cite{BD} and
\cite{Lot1}.

We may now look at the limits of the group actions, which give
\begin{align*}
&  \lim_{s\rightarrow\infty}\phi_{s}^{-1}\circ\gamma_{\left(  a,b,c\right)
}\circ\phi_{s}\left(  x,y,z\right) \\
&  =\lim_{s\rightarrow\infty}\left(  e^{-c}x+s^{-1/2}a,e^{c}y+s^{-1/2}%
b,z+c\right) \\
&  =\left(  e^{-c}x+u,e^{c}y+v,z+c\right)
\end{align*}
if for large $s,$ we choose $a$ and $b$ so that
\begin{align*}
\lim_{s\rightarrow\infty}s^{-1/2}a\left(  s\right)   &  =u,\\
\lim_{s\rightarrow\infty}s^{-1/2}b\left(  s\right)   &  =v.
\end{align*}
If we consider a discrete group which gives a compact manifold quotient, we
see that we may get the limits with real $u,v,$ though not so with $c.$ Thus
we get the group of transformations $\gamma_{\left(  u,v,c\right)  }$ where
$u,v\in\mathbb{R}$ and $c\in\mathbb{Z}$. Clearly this limit is a groupoid
whose orbit space is a circle.

\subsubsection{Cross Curvature Flow}

From \cite{CNS}, we see that the negative cross curvature flow on
$\operatorname{Sol}$ has solutions%
\begin{align*}
B  &  \sim2\sqrt{T_{0}-t},\\
A,C  &  \sim\frac{E_{1}}{\sqrt{T_{0}-t}},\\
A-C  &  \sim E_{2}\sqrt{T_{0}-t}%
\end{align*}
for some positive constants $E_{1}$ and $E_{2}$ and singular time $T_{0}>1.$
Note that all of the curvatures blow up like $\frac{E}{\sqrt{T_{0}-t}},$ so
the solution is Type I. We may now consider the Type I renormalization,%
\begin{align*}
&  g_{T_{0}}\left(  t\right) \\
&  =\lim_{s\rightarrow0}s^{-1/2}\phi_{s}^{\ast}g\left(  T_{0}-st\right) \\
&  =\lim_{s\rightarrow0}s^{-1/2}\left[  \frac{E_{1}}{\sqrt{st}}\left(
e^{z}s^{1/2}dx+e^{-z}s^{1/2}dy\right)  ^{2}+2\sqrt{st}dz^{2}+\frac{E_{1}%
}{\sqrt{st}}\left(  e^{z}s^{1/2}dx-e^{-z}s^{1/2}dy\right)  ^{2}\right] \\
&  =\frac{E_{1}}{\sqrt{t}}\left(  e^{z}dx+e^{-z}dy\right)  ^{2}+2\sqrt
{t}dz^{2}+\frac{E_{1}}{\sqrt{t}}\left(  e^{z}dx-e^{-z}dy\right)  ^{2}\\
&  =\frac{E_{1}}{\sqrt{t}}\left(  e^{2z}dx^{2}+e^{-2z}dy^{2}\right)
+2\sqrt{t}dz^{2},
\end{align*}
where
\[
\phi_{s}\left(  x,y,z\right)  =\left(  s^{1/2}x,s^{1/2}y,z\right)  .
\]
Note that this equals
\[
t^{1/2}\psi_{t}^{\ast}g\left(  1\right)
\]
where
\[
\psi_{t}\left(  x,y,z\right)  =\left(  t^{-1/2}x,t^{-1/2}y,z\right)  .
\]

We may now look at the limit of the group actions,
\begin{align*}
&  \lim_{s\rightarrow0}\phi_{s}^{-1}\circ\gamma_{\left(  a,b,c\right)  }%
\circ\phi_{s}\left(  x,y,z\right) \\
&  =\lim_{s\rightarrow0}\left(  e^{-c}x+s^{-1/2}a,e^{c}y+s^{-1/2}b,z+c\right)
\\
&  =\left(  e^{-c}x+u,e^{c}y+v,z+c\right)
\end{align*}
if
\begin{align*}
\lim_{s\rightarrow0}s^{-1/2}a\left(  s\right)   &  =u,\\
\lim_{s\rightarrow0}s^{-1/2}b\left(  s\right)   &  =v.
\end{align*}
If we consider the limit of a groupoid representing a compact manifold which
is the quotient of $\operatorname{Sol}$ by isometries, then we see that we
cannot choose $a\left(  s\right)  $ and $b\left(  s\right)  $ to be
arbitrarily small (because the group must act properly discontinuously), so we
must have $u=v=0.$ This corresponds to the fact that as we go to the limit,
the arrows of the groupoid send elements of the ball of radius $R$ outside the
ball of radius $2R$ for any $R,$ and so these do not survive in the limit
groupoid. Thus the orbit space looks like a noncompact $\operatorname{Sol}$
with a quotient by the group of transformations $\gamma_{\left(  0,0,c\right)
}$ where $c\in\mathbb{Z}$.

\subsection{$\widetilde{\operatorname{SL}}\left(  2,\mathbb{R}\right)  $}

\begin{figure}
[ptb]
\begin{center}
\fbox{\includegraphics[
natheight=3.000000in,
natwidth=4.499600in,
height=3in,
width=4.4996in
]%
{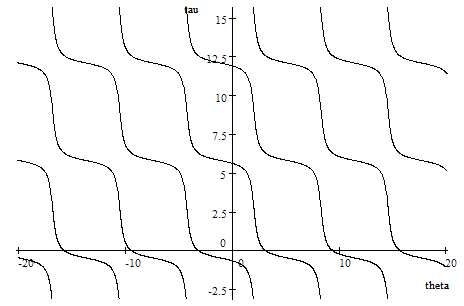}%
}\caption{Solutions of $\cot\frac{\tau}{2}=y\tan\frac{\theta}{2}+x$ for $x=-3$
and $y=2.$}%
\label{wavy lines}%
\end{center}
\end{figure}

The homogeneous geometry $\widetilde{\operatorname{SL}}\left(  2,\mathbb{R}%
\right)  $ is diffeomorphic to $\mathbb{R}^{3},$ but in a nontrivial way. We
construct a coordinate patch on $\widetilde{\operatorname{SL}}\left(
2,\mathbb{R}\right)  $ following the description in \cite{Sco}.
$\operatorname{PSL}\left(  2,\mathbb{R}\right)  $ is the isometry group of the
hyperbolic plane $\mathbb{H}^{2}.$ Thus we see that $\operatorname{PSL}\left(
2,\mathbb{R}\right)  $ acts freely and transitively on $U\mathbb{H}^{2},$ the
unit tangent bundle of $\mathbb{H}^{2},$ and we thus have a diffeomorphism
$\operatorname{PSL}\left(  2,\mathbb{R}\right)  \cong U\mathbb{H}^{2}.$ We
lift this diffeomorphism to the universal cover and easily recognize the
universal cover of the right side as diffeomorphic to $\mathbb{R}^{3}.$ The
lifted diffeomorphism can be described explicitly as follows. We first
describe the map $\Phi:\operatorname{PSL}\left(  2,\mathbb{R}\right)
\rightarrow U\mathbb{H}^{2},$ which is given by the action on the basepoint
$\left(  i,\left[  0\right]  \right)  =\left(  0,1,\left[  0\right]  \right)
\in U\mathbb{H}^{2}$ (where the unit vectors are described by their angles,
$\left[  \theta\right]  =\theta\operatorname{mod}2\pi,$ with respect to the
$x$-axis in the half-plane model):%
\begin{align}
\Phi\left(
\begin{array}
[c]{cc}%
a & b\\
c & d
\end{array}
\right)   &  =\left(  \frac{ai+b}{ci+d},\left[  \tan^{-1}\frac{2cd}%
{d^{2}-c^{2}}\right]  \right)  \nonumber\\
&  =\left(  \frac{ac+bd}{c^{2}+d^{2}},\frac{1}{c^{2}+d^{2}},\left[  \tan
^{-1}\frac{2cd}{d^{2}-c^{2}}\right]  \right)  .\label{PSL multiplication}%
\end{align}
The diffeomorphism is then gotten by lifting the map
\[
\widetilde{\operatorname{SL}}\left(  2,\mathbb{R}\right)  \rightarrow
\operatorname{PSL}\left(  2,\mathbb{R}\right)  \overset{\Phi}{\rightarrow
}U\mathbb{H}^{2}%
\]
to the universal cover of $U\mathbb{H}^{2}.$ For future use, we explicitly
give the inverse of $\Phi$:%
\[
\Phi^{-1}\left(  x,y,\left[  \theta\right]  \right)  =\left[  \frac{1}%
{y^{1/2}}\left(
\begin{array}
[c]{cc}%
x\sin\frac{\theta}{2}+y\cos\frac{\theta}{2} & x\cos\frac{\theta}{2}-y\sin
\frac{\theta}{2}\\
\sin\frac{\theta}{2} & \cos\frac{\theta}{2}%
\end{array}
\right)  \right]  ,
\]
where $\left[  \left(
\begin{array}
[c]{cc}%
\cdot & \cdot\\
\cdot & \cdot
\end{array}
\right)  \right]  $ denotes the equivalence class of the matrix up to
multiplication by $-1.$ We can derive the group multiplication $\left(
a,b,\left[  \tau\right]  \right)  \left(  x,y,\left[  \theta\right]  \right)
$ in $U\mathbb{H}^{2}$ to be
\begin{align}
\left(  a,b,\left[  \tau\right]  \right)  \left(  x,y,\left[  \theta\right]
\right)   &  =\Phi\left(  \Phi^{-1}\left(  a,b,\left[  \tau\right]  \right)
\cdot\Phi^{-1}\left(  x,y,\left[  \theta\right]  \right)  \right)  \nonumber\\
&  =\left(
\begin{array}
[c]{c}%
a+\frac{b\left(  x\cos\tau+\frac{1}{2}\left(  x^{2}+y^{2}-1\right)  \sin
\tau\right)  }{\sin^{2}\frac{\tau}{2}\left[  \left(  x+\cot\frac{\tau}%
{2}\right)  ^{2}+y^{2}\right]  },\\
\frac{by}{\sin^{2}\frac{\tau}{2}\left[  \left(  x+\cot\frac{\tau}{2}\right)
^{2}+y^{2}\right]  },\\
\left[  2\tan^{-1}\left(  \frac{\sin\frac{\theta}{2}\left(  x+\cot\frac{\tau
}{2}\right)  +y\cos\frac{\theta}{2}}{\cos\frac{\theta}{2}\left(  x+\cot
\frac{\tau}{2}\right)  -y\sin\frac{\theta}{2}}\right)  \right]
\end{array}
\right)  \nonumber\\
&  =\left(
\begin{array}
[c]{c}%
a+\frac{b\left(  x\cos\tau+\frac{1}{2}\left(  x^{2}+y^{2}-1\right)  \sin
\tau\right)  }{\sin^{2}\frac{\tau}{2}\left[  \left(  x+\cot\frac{\tau}%
{2}\right)  ^{2}+y^{2}\right]  },\\
\frac{by}{\sin^{2}\frac{\tau}{2}\left[  \left(  x+\cot\frac{\tau}{2}\right)
^{2}+y^{2}\right]  },\\
\left[  \theta+2\tan^{-1}\left(  \frac{y}{x+\cot\frac{\tau}{2}}\right)
\right]
\end{array}
\right)  .\label{multiplication}%
\end{align}
Certainly the universal cover of $U\mathbb{H}^{2}$ is diffeomorphic to
$\mathbb{R}\times\mathbb{R}_{+}\times\mathbb{\mathbb{R}}$, say consisting of
elements like $\left(  x,y,\theta\right)  .$ We may lift the multiplication
map
\[
\mu:\widetilde{U\mathbb{H}}^{2}\times\widetilde{U\mathbb{H}}^{2}\rightarrow
U\mathbb{H}^{2}%
\]
given by (\ref{multiplication}) to a map
\[
\tilde{\mu}:\widetilde{U\mathbb{H}}^{2}\times\widetilde{U\mathbb{H}}%
^{2}\rightarrow\widetilde{U\mathbb{H}}^{2}%
\]
where we specify that $\tilde{\mu}\left(  \left(  0,1,0\right)  ,\left(
0,1,0\right)  \right)  =\left(  0,1,0\right)  .$ The first two coordinates are
unchanged, but we must lift the map to the third coordinates. We denote by
$\tilde{\mu}_{3}$ the lifted map to the third coordinate$.$ We see by
(\ref{multiplication}) that $\tilde{\mu}_{3}=\tilde{\mu}_{3}\left(
\tau,x,y,\theta\right)  .$ To get a handle on the lift, we need to see when%
\[
\cos\frac{\theta}{2}\left(  x+\cot\frac{\tau}{2}\right)  -y\sin\frac{\theta
}{2}=0,
\]
i.e.,
\[
\cot\frac{\tau}{2}=y\tan\frac{\theta}{2}+x.
\]
Note that if $x=0$ and $y=1,$ then the solutions are $\tau=\pi\left(
2k+1\right)  -\theta$ for any $k\in\mathbb{Z}$. For general $x\in\mathbb{R}$
and $y\in\mathbb{R}_{+},$ the solutions are wiggly lines which roughly follow
those lines, as the set of all lines is invariant under translation by
multiples of $2\pi$ in both the up/down and left/right directions (see Figure
\ref{wavy lines}).%
In particular, we see that on these curves $\tau$ is decreasing as a function
of $\theta$ and that the curves are invariant under translations in $\tau$ and
$\theta$ by integer multiples of $2\pi.$ There is one curve which intersects
$\left(  \theta,\tau\right)  =\left(  \pi,0\right)  $ and that curve also
intersects $\left(  0,\pi-2\tan^{-1}x\right)  .$ The behavior can be
understood by looking at blocks $\left(  \theta,\tau\right)  \in\left[
0,2\pi\right]  \times\left[  0,2\pi\right]  $ (in addition to translations of
this curve by multiples of $2\pi$ in both directions). We see that if $\theta$
and $\tau$ are positive, then $\tilde{\mu}_{3}$ is essentially%
\begin{align}
\tilde{\mu}_{3}\left(  \tau,x,y,\theta\right)   &  \approx2\tan^{-1}\left(
\frac{\sin\frac{\theta}{2}\left(  x+\cot\frac{\tau}{2}\right)  +y\cos
\frac{\theta}{2}}{\cos\frac{\theta}{2}\left(  x+\cot\frac{\tau}{2}\right)
-y\sin\frac{\theta}{2}}\right)  +2\pi\left\lfloor \frac{\theta+\tau}{2\pi
}\right\rfloor \label{formula for mu3}\\
&  \approx\theta+2\tan^{-1}\left(  \frac{y}{\left(  x+\cot\frac{\tau}%
{2}\right)  }\right)  +2\pi\left\lfloor \frac{\tau}{2\pi}\right\rfloor
\nonumber
\end{align}
where $\left\lfloor \cdot\right\rfloor $ is the floor (greatest integer less than).

\begin{remark}
Another way to get coordinates on $\widetilde{\operatorname{SL}}\left(
2,\mathbb{R}\right)  $ is to first consider the Iwasawa decomposition of
$\operatorname{SL}\left(  2,\mathbb{R}\right)  ,$ which says that every matrix
in $\operatorname{SL}\left(  2,\mathbb{R}\right)  $ can be written as
\[
\left(
\begin{array}
[c]{cc}%
\cos\theta & \sin\theta\\
-\sin\theta & \cos\theta
\end{array}
\right)  \left(
\begin{array}
[c]{cc}%
y & x\\
0 & \frac{1}{y}%
\end{array}
\right)
\]
where $y>0.$ In order to write down the multiplication, one needs to rewrite
the product in this form again, and then lift to the universal cover. It turns
out that to perform these two operations, the specifics of how $x,$ $y,$ and
$\theta$ change in this coordinate chart is much the same as the way they
behave in the first coordinate chart we gave. This is because the Iwasawa
decomposition of $\operatorname{SL}\left(  2,\mathbb{R}\right)  $ acts in an
easy way on $U\mathbb{H}^{2},$ i.e., the rotation $\left(
\begin{array}
[c]{cc}%
\cos\theta & \sin\theta\\
-\sin\theta & \cos\theta
\end{array}
\right)  $ fixes $i$ and changes the direction of the vector, and the matrix
$\left(
\begin{array}
[c]{cc}%
y & x\\
0 & \frac{1}{y}%
\end{array}
\right)  $ acts by moving $i$ to another point in $\mathbb{H}^{2}$ but fixing
the direction.
\end{remark}

There is a basis of left invariant vector fields given by
\begin{align*}
F_{1}  &  =-\frac{\partial}{\partial\theta},\\
F_{2}  &  =y\cos\theta\frac{\partial}{\partial x}-y\sin\theta\frac{\partial
}{\partial y}+\cos\theta\frac{\partial}{\partial\theta},\\
F_{3}  &  =y\sin\theta\frac{\partial}{\partial x}+y\cos\theta\frac{\partial
}{\partial y}+\sin\theta\frac{\partial}{\partial\theta},
\end{align*}
and we can see easily that
\begin{align*}
\left[  F_{1},F_{2}\right]   &  =F_{3},\\
\left[  F_{2},F_{3}\right]   &  =-F_{1},\\
\left[  F_{3},F_{1}\right]   &  =F_{2}.
\end{align*}
The following is a family of left invariant metrics:
\begin{equation}
A\left(  d\theta-\frac{1}{y}dx\right)  ^{2}+B\left(  \frac{1}{y}\cos\theta
dx-\frac{1}{y}\sin\theta dy\right)  ^{2}+C\left(  \frac{1}{y}\sin\theta
dx+\frac{1}{y}\cos\theta dy\right)  ^{2}. \label{SL(2,R) metrics}%
\end{equation}
If $B=C,$ then the metric is
\[
g=A\left(  d\theta-\frac{1}{y}dx\right)  ^{2}+\frac{B}{y^{2}}\left(
dx^{2}+dy^{2}\right)
\]
and the metric is a bundle over $\mathbb{H}^{2}.$ Note that the isometry group
for the general metrics (\ref{SL(2,R) metrics}) is $\widetilde
{\operatorname{SL}}\left(  2,\mathbb{R}\right)  ,$ but if $B=C,$ there is an
additional part which includes arbitrary translations of $\theta,$ so that the
isometry group looks like $\widetilde{\operatorname{SL}}\left(  2,\mathbb{R}%
\right)  \times\mathbb{R}.$ (Actually, when $B=C$ this is only the identity
component; there are two components due the the isometry $\left(
x,y,\theta\right)  \rightarrow\left(  -x,y,-\theta\right)  .$ For more, see
\cite{Sco}.)

The sectional curvatures are:%
\begin{align*}
K\left(  F_{2}\wedge F_{3}\right)   &  =\frac{-3A^{2}+B^{2}+C^{2}%
-2AB-2BC-2AC}{4ABC},\\
K\left(  F_{3}\wedge F_{1}\right)   &  =\frac{A^{2}-3B^{2}+C^{2}%
-2AB+2BC+2AC}{4ABC},\\
K\left(  F_{1}\wedge F_{2}\right)   &  =\frac{A^{2}+B^{2}-3C^{2}%
+2AB+2BC-2AC}{4ABC}.
\end{align*}
If $B=C,$ the sectional curvatures are:%
\begin{align*}
K\left(  F_{2}\wedge F_{3}\right)   &  =\frac{-3A-4B}{4B^{2}},\\
K\left(  F_{3}\wedge F_{1}\right)   &  =\frac{A}{4B^{2}},\\
K\left(  F_{1}\wedge F_{2}\right)   &  =\frac{A}{4B^{2}}.
\end{align*}

Note that we have taken $F_{1},F_{2},F_{3}$ which are one half that in
\cite{KM} and \cite{CNS}, so out $A,B,C$ are all one fourth the corresponding
coefficients in those papers.

\subsubsection{Ricci Flow}

Under Ricci flow, from \cite{KM} we have that $A$ goes to a constant, $B$ and
$C$ are like $2t$ and
\[
\left\vert B-C\right\vert \leq E_{1}e^{-E_{2}t},
\]
for positive constants $E_{1}$ and $E_{2}$. We see that the sectional
curvatures can be written as%
\begin{align*}
K\left(  F_{2}\wedge F_{3}\right)   &  =\frac{\left(  B-C\right)
^{2}-A\left(  3A+2B+2C\right)  }{4ABC},\\
K\left(  F_{3}\wedge F_{1}\right)   &  =\frac{\left(  A-\left(  B-C\right)
\right)  ^{2}-4B\left(  B-C\right)  }{4ABC},\\
K\left(  F_{1}\wedge F_{2}\right)   &  =\frac{\left(  A+\left(  B-C\right)
\right)  ^{2}+4C\left(  B-C\right)  }{4ABC},
\end{align*}
and so
\begin{align*}
K\left(  F_{2}\wedge F_{3}\right)   &  \sim\frac{E_{3}}{t},\\
K\left(  F_{3}\wedge F_{1}\right)   &  \sim\frac{E_{4}}{t^{2}},\\
K\left(  F_{1}\wedge F_{2}\right)   &  \sim\frac{E_{5}}{t^{2}}%
\end{align*}
for large $t,$ where $E_{3},E_{4},E_{5}$ are constants. Thus the solution is
Type III.

We may do a Type III rescaling and pull back by the diffeomorphism
\[
\phi_{s}\left(  x,y,\theta\right)  =\phi_{s}\left(  x,y,s^{1/2}\theta\right)
\]
to get the limit
\[
\lim_{s\rightarrow\infty}\frac{1}{s}\phi_{s}^{\ast}g\left(  st\right)
=d\theta^{2}+\frac{2t}{y^{2}}\left(  dx^{2}+dy^{2}\right)  ,
\]
which is an expanding soliton on $\mathbb{H}^{2}\times\mathbb{R}$.

Let's look at the evolution of the action of the isometry group. We see that
\[
\phi_{s}^{-1}\left(  a,b,\tau\right)  \phi_{s}\left(  x,y,\theta\right)
=\left(
\begin{array}
[c]{c}%
a+\frac{b\left(  x\cos\tau+\frac{1}{2}\left(  x^{2}+y^{2}-1\right)  \sin
\tau\right)  }{\sin^{2}\frac{\tau}{2}\left[  \left(  x+\cot\frac{\tau}%
{2}\right)  ^{2}+y^{2}\right]  },\\
\frac{by}{\sin^{2}\frac{\tau}{2}\left[  \left(  x+\cot\frac{\tau}{2}\right)
^{2}+y^{2}\right]  },\\
s^{-1/2}\tilde{\mu}_{3}\left(  \tau,x,y,s^{1/2}\theta\right)
\end{array}
\right)  .
\]
We note that for large $s,$
\[
s^{-1/2}\tilde{\mu}_{3}\left(  \tau,x,y,s^{1/2}\theta\right)  \approx
\theta+s^{-1/2}\tan^{-1}\left(  \frac{y}{x+\cot\frac{\tau}{2}}\right)  +2\pi
s^{-1/2}\left\lfloor \frac{\tau}{2\pi}\right\rfloor .
\]
Thus if $\tau\left(  s\right)  $ is bounded, then as $s\rightarrow\infty$ we
get elements that look like:%
\[
\gamma_{\left(  a,b,\tau\right)  }\left(  x,y,\theta\right)  =\left(
\begin{array}
[c]{c}%
a+\frac{b\left(  x\cos\tau+\frac{1}{2}\left(  x^{2}+y^{2}-1\right)  \sin
\tau\right)  }{\sin^{2}\frac{\tau}{2}\left[  \left(  x+\cot\frac{\tau}%
{2}\right)  ^{2}+y^{2}\right]  },\\
\frac{by}{\sin^{2}\frac{\tau}{2}\left[  \left(  x+\cot\frac{\tau}{2}\right)
^{2}+y^{2}\right]  },\\
\theta
\end{array}
\right)  .
\]
If $\tau\left(  s\right)  \rightarrow\infty,$ then we may have other elements.
For the first two components to make sense, we need to take a sequence where
$\tau_{i}=\tau_{0}+2\pi k_{i}$ for $k_{i}\in\mathbb{Z}$. If we take such a
sequence of $\tau_{i}$ and a sequence of $s_{i}\rightarrow\infty$ such that
\[
\lim_{i\rightarrow\infty}2\pi s_{i}^{-1/2}\left\lfloor \frac{\tau_{i}}{2\pi
}\right\rfloor =u
\]
for $u\in\mathbb{R}$, we see that we can get as a limit elements that look
like
\[
\gamma_{\left(  a,b,\tau,u\right)  }\left(  x,y,\theta\right)  =\left(
\begin{array}
[c]{c}%
a+\frac{b\left(  x\cos\tau+\frac{1}{2}\left(  x^{2}+y^{2}-1\right)  \sin
\tau\right)  }{\sin^{2}\frac{\tau}{2}\left[  \left(  x+\cot\frac{\tau}%
{2}\right)  ^{2}+y^{2}\right]  },\\
\frac{by}{\sin^{2}\frac{\tau}{2}\left[  \left(  x+\cot\frac{\tau}{2}\right)
^{2}+y^{2}\right]  },\\
\theta+u
\end{array}
\right)  .
\]
In particular, if $\tau_{0}=0$ and $a=0,b=1$ then we get the translations
\[
\gamma_{\left(  0,1,0,u\right)  }\left(  x,y,\theta\right)  =\left(
\begin{array}
[c]{c}%
x,\\
y,\\
\theta+u
\end{array}
\right)
\]
for any $u\in\mathbb{R}$. Note that for simplicity of the formula in
(\ref{formula for mu3}) we assumed that $\tau$ and $\theta$ are positive, but
one can also do the general case with a careful analysis of the lifted map
$\tilde{\mu}_{3}.$

Thus the isometry group converges the standard action of $\operatorname{PSL}%
\left(  2,R\right)  $ on $\mathbb{H}^{2}$ in the first two coordinates and a
continuous action in the last coordinate, i.e., to $\operatorname{PSL}\left(
2,\mathbb{R}\right)  \times\mathbb{R}$.

We may consider a compact quotient $\widetilde{\operatorname{SL}}\left(
2,\mathbb{R}\right)  /\Gamma$ where $\Gamma$ acts properly discontinuously
(and freely if the quotient is a manifold). A properly discontinuous group
action may still contain the $\mathbb{R}$ action since $\tau$ may grow without
bound, and so the limit group will still be continuous. Thus the orbit space
of the groupoid will be a two-dimensional quotient of $\mathbb{H}^{2}%
\times\mathbb{R}$.

\subsubsection{Cross Curvature Flow}

The negative cross curvature flow exhibits two different types of behavior, so
we divide it into two cases.

\paragraph{Case 1: $B_{0}=C_{0}.$}

In this case, the solution exists for all $t\in\lbrack0,\infty)$ and $B\left(
t\right)  =C\left(  t\right)  $ for all $t.$ The solutions are%
\[
A\sim A_{\infty},\;\;\;B=C\sim\left(  \frac{3}{2}A_{\infty}t\right)  ^{1/3},
\]
where $A_{\infty}>0$ is the limit of $A,$ which decreases monotonically to it.
The sectional curvatures are%
\begin{align*}
K\left(  F_{2}\wedge F_{3}\right)   &  =\frac{-3A-4B}{4B^{2}}\sim-\frac
{1}{\left(  \frac{3}{2}A_{\infty}t\right)  ^{1/3}},\\
K\left(  F_{3}\wedge F_{1}\right)   &  =\frac{A}{4B^{2}}\sim\frac{A_{\infty}%
}{4\left(  \frac{3}{2}A_{\infty}t\right)  ^{2/3}},\\
K\left(  F_{1}\wedge F_{2}\right)   &  =\frac{A}{4B^{2}}\sim\frac{A_{\infty}%
}{4\left(  \frac{3}{2}A_{\infty}t\right)  ^{2/3}}.
\end{align*}
We see that this is a Type IIb solution since one sectional curvature does not
decay faster than $t^{-1/2}.$ We may take the geometric limit%
\begin{align*}
&  \lim_{s\rightarrow\infty}s^{-1/3}\phi_{s}^{\ast}g\left(  s^{2/3}\left(
t-1\right)  +s\right) \\
&  =\lim_{s\rightarrow\infty}\left[  s^{-1/3}A_{\infty}\left(  s^{1/6}%
d\theta-\frac{1}{y}dx\right)  ^{2}+s^{-1/3}\left(  \frac{3}{2}A_{\infty
}\right)  ^{1/3}\left(  s^{2/3}\left(  t-1\right)  +s\right)  ^{1/3}\left(
\frac{1}{y^{2}}dx^{2}+\frac{1}{y^{2}}dy^{2}\right)  \right] \\
&  =A_{\infty}d\theta^{2}+\left(  \frac{3}{2}A_{\infty}\right)  ^{1/3}\frac
{1}{y^{2}}\left(  dx^{2}+dy^{2}\right)
\end{align*}
where
\[
\phi_{s}\left(  x,y,\theta\right)  =\left(  x,y,s^{1/6}\theta\right)  .
\]
This is the solution $\mathbb{H}^{2}\times\mathbb{R}$, which is a fixed point
of the flow. (In fact, since the cross curvature tensor consists of products
of two of the sectional curvatures, any product metric is a fixed point.) The
compact quotients behave the same way as the Ricci flow case.

\paragraph{Case 2: $B_{0}>C_{0}$}

In this case, \cite{CNS} shows that
\[
A,B\sim E\left(  T_{0}-t\right)  ^{-1/2},\;\;C\sim2\sqrt{T_{0}-t},
\]
for some singularity time $T_{0}>1$ and that
\begin{align*}
\lim_{t\rightarrow T_{0}}C~K\left(  F_{2}\wedge F_{3}\right)   &  =-1,\\
\lim_{t\rightarrow T_{0}}C~K\left(  F_{3}\wedge F_{1}\right)   &  =-1,\\
\lim_{t\rightarrow T_{0}}C~K\left(  F_{1}\wedge F_{3}\right)   &  =1.
\end{align*}
Thus the sectional curvatures blow up like
\[
\frac{1}{2\sqrt{T_{0}-t}},
\]
and the singularity is Type I. We look at the Type I limit:%

\begin{align*}
&  \lim_{s\rightarrow0}s^{-1/2}\phi_{s}^{\ast}g\left(  T_{0}-st\right) \\
&  =\lim_{s\rightarrow0}\left[
\begin{array}
[c]{r}%
Es^{-1}t^{-1/2}\left(  s^{1/2}d\theta-\frac{1}{y}s^{1/2}dx\right)
^{2}+Es^{-1}t^{-1/2}\left(  \frac{1}{y}\cos\left(  s^{1/2}\theta\right)
s^{1/2}dx-\frac{1}{y}\sin\left(  s^{1/2}\theta\right)  dy\right)  ^{2}\\
+2t^{1/2}\left(  \frac{1}{y}\sin\left(  s^{1/2}\theta\right)  s^{1/2}%
dx+\frac{1}{y}\cos\left(  s^{1/2}\theta\right)  dy\right)  ^{2}%
\end{array}
\right] \\
&  =Et^{-1/2}\left(  d\theta-\frac{1}{y}dx\right)  ^{2}+Et^{-1/2}\left(
\frac{1}{y}dx-\frac{\theta}{y}dy\right)  ^{2}+2t^{1/2}\left(  \frac{1}%
{y}dy\right)  ^{2},
\end{align*}
where
\[
\phi_{s}\left(  x,y,\theta\right)  =\left(  s^{1/2}x,y,s^{1/2}\theta\right)
.
\]
We claim that this is actually the XC soliton on Sol. We can see this if we
pull back by%
\[
\psi\left(  x,y,\theta\right)  =\left(  e^{y}x,e^{y},\theta+x\right)  =\left(
\tilde{x},\tilde{y},\tilde{\theta}\right)
\]
to get
\begin{align*}
&  Et^{-1/2}\left(  d\tilde{\theta}-\frac{1}{\tilde{y}}d\tilde{x}\right)
^{2}+Et^{-1/2}\left(  \frac{1}{\tilde{y}}d\tilde{x}-\frac{\tilde{\theta}%
}{\tilde{y}}d\tilde{y}\right)  ^{2}+2t^{1/2}\left(  \frac{1}{\tilde{y}}%
d\tilde{y}\right)  ^{2}\\
&  =Et^{-1/2}\left(  d\theta-xdy\right)  ^{2}+Et^{-1/2}\left(  dx-\theta
dy\right)  ^{2}+2t^{1/2}dy^{2},
\end{align*}
which can easily be transformed by rescaling coordinates to the metric
(\ref{alt sol metric}).

We let
\begin{align*}
\psi_{s}\left(  x,y,\theta\right)   &  =\phi_{s}\circ\psi\left(
x,y,\theta\right)  =\left(  s^{1/2}e^{y}x,e^{y},s^{1/2}\left(  \theta
+x\right)  \right) \\
\psi_{s}^{-1}\left(  x,y,\theta\right)   &  =\left(  s^{-1/2}\frac{x}{y},\log
y,s^{-1/2}\left(  \theta-\frac{x}{y}\right)  \right)
\end{align*}
Now consider what happens to the group action in the limit. We compute%
\begin{align*}
&  \psi_{s}^{-1}\left(  a,b,\tau\right)  \psi_{s}\left(  x,y,\theta\right) \\
&  =\psi_{s}^{-1}\left(
\begin{array}
[c]{c}%
a+\frac{b\left(  s^{1/2}e^{y}x\cos\tau+\frac{1}{2}\left(  se^{2y}x^{2}%
+e^{2y}-1\right)  \sin\tau\right)  }{\sin^{2}\frac{\tau}{2}\left[  \left(
s^{1/2}e^{y}x+\cot\frac{\tau}{2}\right)  ^{2}+y^{2}\right]  },\\
\frac{be^{y}}{\sin^{2}\frac{\tau}{2}\left[  \left(  s^{1/2}e^{y}x+\cot
\frac{\tau}{2}\right)  ^{2}+e^{2y}\right]  },\\
\tilde{\mu}_{3}\left(  \tau,s^{1/2}e^{y}x,e^{y},s^{1/2}\left(  \theta
+x\right)  \right)
\end{array}
\right) \\
&  =\left(
\begin{array}
[c]{c}%
\alpha,\\
\beta,\\
\gamma,
\end{array}
\right)
\end{align*}
where%
\begin{align*}
\alpha_{s}\left(  a,b,\tau,x,y\right)   &  =s^{-1/2}\frac{a}{b}e^{y}\left(
\left(  s^{1/2}x\sin\frac{\tau}{2}+e^{-y}\cos\frac{\tau}{2}\right)  ^{2}%
+\sin^{2}\frac{\tau}{2}\right) \\
&  +\left(  x\cos\tau+\frac{1}{2}\left(  s^{1/2}e^{y}x^{2}+s^{-1/2}%
e^{y}-s^{-1/2}e^{-y}\right)  \sin\tau\right)  ,
\end{align*}%
\[
\beta_{s}\left(  b,\tau,x,y\right)  =y+\log b-\log\left(  \left(  s^{1/2}%
e^{y}x\sin\frac{\tau}{2}+\cos\frac{\tau}{2}\right)  ^{2}+e^{2y}\sin^{2}%
\frac{\tau}{2}\right)  ,
\]
and%
\[
\gamma_{s}\left(  a,b,\tau,x,y,\theta\right)  =s^{-1/2}\tilde{\mu}_{3}\left(
\tau,s^{1/2}e^{y}x,e^{y},s^{1/2}\left(  \theta+x\right)  \right)  -\alpha
_{s}\left(  a,b,\tau,x,y\right)  .
\]
We see immediately that since $\alpha_{s}$ cannot become unbounded as
$s\rightarrow0,$ for $\gamma_{s}$ to not be unbounded, we need
\[
s^{-1/2}\tilde{\mu}_{3}\left(  \tau,s^{1/2}e^{y}x,e^{y},s^{1/2}\left(
\theta+x\right)  \right)
\]
to stay bounded. Since the last term for $\tilde{\mu}_{3}$ goes to zero, this
is a restriction on our choice of $\tau\left(  s\right)  .$ In particular, we
need that $\tau$ stays relatively close to zero. Knowing this, we can look
more precisely at the formula for this term, which, for positive $\tau$, is%
\begin{align*}
&  s^{-1/2}\tilde{\mu}_{3}\left(  \tau,s^{1/2}e^{y}x,e^{y},s^{1/2}\left(
\theta+x\right)  \right) \\
&  \approx\theta+x+2s^{-1/2}\tan^{-1}\left(  \frac{e^{y}\tan\frac{\tau}{2}%
}{\left(  s^{1/2}e^{y}x\tan\frac{\tau}{2}+1\right)  }\right)  +2\pi
s^{-1/2}\left\lfloor \frac{\tau}{2\pi}\right\rfloor .
\end{align*}
Since $\tau$ is close to zero, the last term is always zero, and for the
second to last term to not go to infinity, we need
\[
\lim_{s\rightarrow0}s^{-1/2}\tau=2u,
\]
for some $u\in\mathbb{R}$, in which case we get
\[
\lim_{s\rightarrow0}\gamma_{s}\left(  a,b,\tau,x,y,\theta\right)
=\theta+x+2e^{y}u-\lim_{s\rightarrow0}\alpha_{s}\left(  a,b,\tau,x,y\right)
.
\]
Considering $\beta_{s}$ and the fact that $\tau\left(  s\right)  \rightarrow0$
as $s\rightarrow0,$ we see that $b$ must converge to a positive number as
$s\rightarrow0,$ so we might as well assume $b\left(  s\right)  \rightarrow
e^{d}.$ Thus
\[
\lim_{s\rightarrow0}\beta_{s}\left(  b,\tau,x,y\right)  =y+d.
\]
Finally, if we have $v\in\mathbb{R}$ such that%
\[
\lim_{s\rightarrow0}s^{-1/2}\frac{a\left(  s\right)  }{b\left(  s\right)
}=u+v,
\]
we have%
\begin{align*}
\lim_{s\rightarrow0}\alpha_{s}\left(  a,b,\tau,x,y\right)   &  =x+e^{-y}%
v+\left(  e^{y}-e^{-y}\right)  u\\
&  =x+e^{y}u+e^{-y}v.
\end{align*}
Thus,
\[
\lim_{s\rightarrow0}\gamma_{s}\left(  a,b,\tau,x,y,\theta\right)
=\theta+e^{y}u-e^{-y}v.
\]
The limit group actions consist of maps
\[
\gamma_{\left(  u,v,d\right)  }\left(  x,y,\theta\right)  =\left(
x+e^{y}u+e^{-y}v,y+d,\theta+e^{y}u-e^{-y}v\right)  .
\]
Note that this has the form of the isometries described by
(\ref{sol isometries}) modulo the change of coordinates defined by scaling
$y.$

Now, if we begin with a groupoid representing a compact manifold quotient of
$\widetilde{\operatorname{SL}}\left(  2,\mathbb{R}\right)  $, the group which
determines the arrows in the groupoid must act properly discontinuously. We
see that this implies that $u$ and $v$ are zero, and $d$ takes discrete
values. Thus the limit is a noncompact quotient of $\operatorname{Sol}$ with
no collapsing.

\subsection{$\widetilde{\operatorname{Isom}}\left(  \mathbb{E}^{2}\right)  $}

The group $\operatorname{Isom}\left(  \mathbb{E}^{2}\right)  $ consists of
group elements
\[
\left(  x^{1},x^{2}\right)  \rightarrow A_{\theta}\left(
\begin{array}
[c]{c}%
x^{1}\\
x^{2}%
\end{array}
\right)  +\left(
\begin{array}
[c]{c}%
x\\
y
\end{array}
\right)  ,
\]
where $A_{\theta}$ is a rotation by angle $\theta$. Thus the universal cover
is diffeomorphic to $\mathbb{R}^{3}$ and has coordinates $\left(
x,y,\theta\right)  .$ The group multiplication is
\[
\left(  a,b,\tau\right)  \left(  x,y,\theta\right)  =\left(  x\cos\tau
+y\sin\tau+a,-x\sin\tau+y\cos\tau+b,\theta+\tau\right)  .
\]
This group has a left invariant frame%
\[
F_{1}=\sin\theta\frac{\partial}{\partial x}+\cos\theta\frac{\partial}{\partial
y},\qquad F_{2}=\cos\theta\frac{\partial}{\partial x}-\sin\theta\frac
{\partial}{\partial y},\qquad F_{3}=\frac{\partial}{\partial\theta}%
\]
with%
\begin{align*}
\left[  F_{2},F_{3}\right]   &  =F_{1}\\
\left[  F_{3},F_{1}\right]   &  =F_{2}%
\end{align*}
the only nonzero brackets. Thus the following are left invariant metrics:
\[
g=A\left(  \sin\theta dx+\cos\theta dy\right)  ^{2}+B\left(  \cos\theta
dx-\sin\theta dy\right)  ^{2}+Cd\theta^{2}.
\]
We note that if $A=B,$ then the metric is Euclidean. It is clear by changing
coordinates by scaling $x$ and $y$ that this is really a two parameter family
of metrics up to diffeomorphism.

The sectional curvatures for these metrics are%
\begin{align*}
K\left(  F_{2}\wedge F_{3}\right)   &  =\frac{\left(  B-A\right)  \left(
B+3A\right)  }{4ABC},\\
K\left(  F_{3}\wedge F_{1}\right)   &  =\frac{\left(  A-B\right)  \left(
A+3B\right)  }{4ABC},\\
K\left(  F_{1}\wedge F_{2}\right)   &  =\frac{\left(  A-B\right)  ^{2}}{4ABC}.
\end{align*}

Note that this $F_{3}$ is one half that used in \cite{KM} and \cite{CNS}, so
our $C$ is 1/4 the corresponding coefficient in those papers.

\subsubsection{Ricci Flow}

From \cite{KM}, we see that the solution to the Ricci flow looks like%
\begin{align*}
A,B  &  \sim E_{1},\\
C  &  \sim E_{2},
\end{align*}
where $E_{1}=\sqrt{A_{0}B_{0}}$ and $E_{2}=\frac{C_{0}}{2}\left(  \sqrt
{\frac{A_{0}}{B_{0}}}+\sqrt{\frac{B_{0}}{A_{0}}}\right)  .$ From the work in
\cite{KM}, we easily see that%
\[
\frac{d}{dt}\left(  A-B\right)  =-\left(  A-B\right)  \frac{\left(
A+B\right)  }{C}\left(  \frac{1}{B}+\frac{1}{A}\right)  ,
\]
and so
\[
A-B\sim E_{4}e^{-E_{3}t},
\]
where%
\begin{align*}
E_{3}  &  =\frac{4}{E_{2}},\\
E_{4}  &  =\left(  A_{0}-B_{0}\right)  .
\end{align*}
The sectional curvatures are
\begin{align*}
K\left(  F_{2}\wedge F_{3}\right)   &  \sim-\frac{4E_{4}E_{1}}{E_{1}^{2}E_{2}%
}e^{-E_{3}t},\\
K\left(  F_{3}\wedge F_{1}\right)   &  \sim\frac{4E_{4}E_{1}}{E_{1}^{2}E_{2}%
}e^{-E_{3}t},\\
K\left(  F_{1}\wedge F_{2}\right)   &  \sim\frac{E_{4}^{2}}{E_{1}^{2}E_{2}%
}e^{-2E_{3}t}.
\end{align*}
This is a Type III solution. It is clear that the Type III limit is Euclidean
space since%
\[
g=A\left(  dx^{2}+dy^{2}\right)  +\left(  B-A\right)  \left(  \cos\theta
dx-\sin\theta dy\right)  ^{2}+Cd\theta^{2}%
\]
and so
\begin{align*}
&  g_{\infty}\left(  t\right) \\
&  =\lim_{s\rightarrow\infty}s^{-1}\phi_{s}^{\ast}g\left(  st\right) \\
&  =\lim_{s\rightarrow\infty}s^{-1}\left[  E_{1}s\left(  dx^{2}+dy^{2}\right)
+E_{4}e^{-E_{3}st}s\left(  \cos s^{1/2}\theta~dx-\sin s^{1/2}\theta~dy\right)
^{2}+E_{2}sd\theta^{2}\right] \\
&  =E_{1}\left(  dx^{2}+dy^{2}\right)  +E_{2}d\theta^{2},
\end{align*}
where
\[
\phi_{s}\left(  x,y,\theta\right)  =\left(  s^{1/2}x,s^{1/2}y,s^{1/2}%
\theta\right)  .
\]

One might try to construct a different geometric limit by choosing a different
rescaling, for instance the following:%
\begin{align*}
g_{s}\left(  t\right)   &  =e^{-E_{3}s}\psi_{s}^{\ast}g\left(  e^{E_{3}%
s}\left(  t-1\right)  +s\right) \\
&  =e^{-E_{3}s}E_{1}e^{E_{3}s}\left(  dx^{2}+dy^{2}\right) \\
&  \;\;\;\;\;+e^{-E_{3}s}E_{4}e^{-E_{3}\left(  e^{E_{3}s}\left(  t-1\right)
+s\right)  }e^{E_{3}s}\left(  \cos\left(  e^{E_{3}s/2}\theta\right)
~dx-\sin\left(  e^{E_{3}s/2}\theta\right)  ~dy\right)  ^{2}\\
&  \;\;\;\;\;+e^{-E_{3}s}E_{2}e^{E_{3}s}d\theta^{2}\\
&  =E_{1}\left(  dx^{2}+dy^{2}\right)  +E_{4}e^{-E_{3}\left(  e^{E_{3}%
s}\left(  t-1\right)  +s\right)  }\left(  \cos\left(  e^{E_{3}s/2}%
\theta\right)  ~dx-\sin\left(  e^{E_{3}s/2}\theta\right)  ~dy\right)  ^{2}\\
&  \;\;\;\;\;+E_{2}d\theta^{2}%
\end{align*}
where
\[
\psi_{s}\left(  x,y,\theta\right)  =\left(  e^{E_{3}s/2}x,e^{E_{3}%
s/2}y,e^{E_{3}s/2}\theta\right)  .
\]
Notice that as $s\rightarrow\infty,$ this also converges to Euclidean space.

Under the Type III limit, the limit of the group actions is%
\begin{align*}
&  \lim_{s\rightarrow\infty}\phi_{s}^{-1}\left[  \gamma_{\left(
a,b,\tau\right)  }\phi_{s}\left(  x,y,\theta\right)  \right] \\
&  =\lim_{s\rightarrow\infty}\left(  x\cos\tau+y\sin\tau+s^{-1/2}a,-x\sin
\tau+y\cos\tau+s^{-1/2}b,\theta+s^{-1/2}\tau\right) \\
&  =\left(  x\cos\tau+y\sin\tau+u,-x\sin\tau+y\cos\tau+v,\theta+w\right)
\end{align*}
if we choose $a\left(  s\right)  $ and $b\left(  s\right)  $ so that
\begin{align*}
\lim_{s\rightarrow\infty}s^{-1/2}a\left(  s\right)   &  =u,\\
\lim_{s\rightarrow\infty}s^{-1/2}b\left(  s\right)   &  =v,\\
\lim_{s\rightarrow\infty}s^{-1/2}\tau\left(  s\right)   &  =w,
\end{align*}
for any $u,v,w\in\mathbb{R}$ and we choose $\tau\left(  s\right)  $ so that it
is growing in multiples of $2\pi$ (so $\cos\tau$ and $\sin\tau$ still make
sense). Thus the limit group is $\operatorname{Isom}\mathbb{E}^{2}%
\times\mathbb{R}$. For a compact quotient, we may still find $a,b,\tau$ that
grow as desired, so we still get the whole group in the limit. Thus the orbit
space of the Type III limit is a point.

\subsubsection{Cross Curvature Flow}

The solution to --XCF is
\begin{align*}
A  &  \sim E_{1},\\
B  &  \sim E_{1},\\
C  &  \sim\frac{2E_{2}}{E_{1}}\sqrt{6}t^{1/3},
\end{align*}
with
\[
A-B\sim E_{2}t^{-1/6}.
\]
Thus the sectional curvatures satisfy
\begin{align*}
K\left(  F_{2}\wedge F_{3}\right)   &  \sim-\frac{1}{2\sqrt{6}t^{1/2}},\\
K\left(  F_{3}\wedge F_{1}\right)   &  \sim\frac{1}{2\sqrt{6}t^{1/2}},\\
K\left(  F_{1}\wedge F_{2}\right)   &  \sim\frac{E_{2}}{8E_{1}\sqrt{6}t^{2/3}%
},
\end{align*}
and the solution is Type III. We may compute the Type III limit of the
rescaled solutions
\begin{align*}
g_{s}\left(  t\right)   &  =s^{-1/2}\phi_{s}^{\ast}g\left(  st\right) \\
&  =s^{-1/2}E_{1}s^{1/2}\left(  dx^{2}+dy^{2}\right) \\
&  \;\;\;\;+s^{-1/2}E_{2}s^{-1/6}t^{-1/6}s^{1/2}\left(  \cos s^{1/12}%
\theta~dx-\sin s^{1/12}\theta~dy\right)  ^{2}\\
&  \;\;\;\;+s^{-1/2}\frac{2E_{2}}{E_{1}}\sqrt{6}s^{1/3}t^{1/3}s^{1/6}%
d\theta^{2},
\end{align*}
where
\[
\phi_{s}\left(  x,y,\theta\right)  =\left(  s^{1/4}x,s^{1/4}y,s^{1/12}%
\theta\right)  .
\]
The limit as $s\rightarrow\infty$ is
\[
g_{\infty}\left(  t\right)  =E_{1}\left(  dx^{2}+dy^{2}\right)  +\frac{2E_{2}%
}{E_{1}}\sqrt{6}t^{1/3}\left(  \frac{1}{2}d\theta\right)  ^{2}.
\]
This is obviously Euclidean space.

The pulled back group actions look like%
\begin{align*}
&  \phi_{s}^{-1}\circ\gamma_{\left(  a,b,\tau\right)  }\circ\phi_{s}\left(
x,y,\theta\right) \\
&  =\left(  x\cos\tau+y\sin\tau+s^{-1/4}a,-x\sin\tau+y\cos\tau+s^{-1/4}%
b,\theta+s^{-1/12}\tau\right)  .
\end{align*}
So in the limit, we may take
\begin{align*}
\lim_{s\rightarrow\infty}s^{-1/4}a\left(  s\right)   &  =u\\
\lim_{s\rightarrow\infty}s^{-1/4}b\left(  s\right)   &  =v\\
\lim_{s\rightarrow\infty}s^{-1/12}\tau\left(  s\right)   &  =w
\end{align*}
to get group actions
\[
\gamma_{\left(  \tau,u,v,w\right)  }\left(  x,y,\theta\right)  =\left(
x\cos\tau+y\sin\tau+u,-x\sin\tau+y\cos\tau+v,\theta+w\right)
\]
if we take the limit so that $\tau\left(  s\right)  $ is growing only by
multiples of $2\pi$ so that $\sin\tau$ and $\cos\tau$ still make sense. Note
that even if the original groupoid comes from a compact quotient, we still get
the entirety of the group since we can let $a,b,\tau$ grow (these elements
exist in the lattice). Thus the limit has an orbit space of a point.

\section{Summary}

We may summarize the results about the limits of compact quotients of
homogeneous geometries as follows. See also the tables in Figures
\ref{RF table} and \ref{XCF table} for a summary. In the tables, the geometry
limit is the limit of the universal covers (or limit of $G^{\left(  0\right)
}$) and the dimension (dim) is the dimension of the orbit space.

\begin{figure}[ptb]%
\begin{tabular}
[c]{|c|c|c|c|c|c|}\hline
Geometry & Soliton & Sing. & Geometry & Collapsing & Compact\\
&  & Type & Limit & limit (dim) & limit\\\hline
$\operatorname{Nil}$ & Yes & III & $\operatorname{Nil}$ & Yes $\left(
0\right)  $ & Yes\\\hline
$\operatorname{Sol}$ & Yes & III & $\operatorname{Sol}$ & Yes $\left(
1\right)  $ & Yes\\\hline
$\left.  \widetilde{\operatorname{SL}}\left(  2,\mathbb{R}\right)  \right.  $
& No & III & $\mathbb{H}^{2}\times\mathbb{R}$ & Yes $\left(  2\right)  $ &
Yes\\\hline
$\left.  \widetilde{\operatorname{Isom}}\left(  \mathbb{E}^{2}\right)
\right.  $ & Yes & III & $\mathbb{E}^{3}$ & Yes $\left(  0\right)  $ &
Yes\\\hline
\end{tabular}
\caption{Summary of limits of Ricci flow}%
\label{RF table}%
\end{figure}

\begin{figure}[ptb]%
\begin{tabular}
[c]{|c|c|c|c|c|c|}\hline
Geometry & Soliton & Sing. & Geometry & Collapsing & Compact\\
&  & Type & Limit & limit (dim) & limit\\\hline
$\operatorname{Nil}$ & Yes & III & $\operatorname{Nil}$ & Yes $\left(
0\right)  $ & Yes\\\hline
$\operatorname{Sol}$ & Yes & I & $\operatorname{Sol}$ & No & No\\\hline
$\mathstrut\left.  \widetilde{\operatorname{SL}}\left(  2,\mathbb{R}\right)
\right.  ,\strut$ $B=C$ & No & IIb & $\mathbb{H}^{2}\times\mathbb{R}$ & Yes
$\left(  2\right)  $ & Yes\\\hline
$\left.  \widetilde{\operatorname{SL}}\left(  2,\mathbb{R}\right)  \right.
,\strut$ $B\neq C$ & No & I & $\operatorname{Sol}$ & No & No\\\hline
$\mathstrut\left.  \widetilde{\operatorname{Isom}}\left(  \mathbb{E}%
^{2}\right)  \right.  \strut$ & Yes & III & $\mathbb{E}^{3}$ & Yes $\left(
0\right)  $ & Yes\\\hline
\end{tabular}
\caption{Summary of limits of negative cross curvature flow}%
\label{XCF table}%
\end{figure}

\begin{theorem}
[Lott \cite{Lot1}]The solutions of Ricci flow on $\operatorname{Nil},$
$\operatorname{Sol},$ $\widetilde{\operatorname{SL}}\left(  2,\mathbb{R}%
\right)  ,$ and $\widetilde{\operatorname{Isom}}\left(  \mathbb{E}^{2}\right)
$ are all of Type III. There are soliton solutions on $\operatorname{Nil},$
$\operatorname{Sol},$ and $\widetilde{\operatorname{Isom}}\left(
\mathbb{E}^{2}\right)  $ (this soliton is $\mathbb{E}^{3},$ three-dimensional
Euclidean space) and the Type III limits converge to these geometries. The
Type III limit of $\widetilde{\operatorname{SL}}\left(  2,\mathbb{R}\right)  $
is $\mathbb{H}^{2}\times\mathbb{R}$. Compact homogeneous manifolds with these
geometries all collapse and stay compact.
\end{theorem}

\begin{theorem}
The solutions of cross curvature flow are as follows:

\begin{itemize}
\item $\operatorname{Nil}$ admits a XC soliton metric. Its solution is Type
III and compact manifolds modeled on $\operatorname{Nil}$ converge to compact,
collapsed quotients of $\operatorname{Nil}$ in the Type III limit.

\item $\operatorname{Sol}$ admits a XC soliton metric. Its solution is Type
IIa and compact manifolds modeled on $\operatorname{Sol}$ converge to
noncollapsed, noncompact quotients of $\operatorname{Sol}.$

\item $\widetilde{\operatorname{SL}}\left(  2,\mathbb{R}\right)  $ does not
seem to admit a soliton metric. If $B=C$ (i.e., if the metric is a Riemannian
submersion over $\mathbb{H}^{2}$), then the singularity is Type IIb and the
geometric limits of compact manifolds modeled on this type of $\widetilde
{\operatorname{SL}}\left(  2,\mathbb{R}\right)  $ are collapsed, compact
quotients of $\mathbb{H}^{2}\times\mathbb{R}$. If $B\neq C,$ then the
singularity is Type I and the type I limits of compact manifolds modeled on
this type of $\widetilde{\operatorname{SL}}\left(  2,\mathbb{R}\right)  $
converge to noncollapsed, noncompact quotients of $\operatorname{Sol}$.

\item $\widetilde{\operatorname{Isom}}\left(  \mathbb{E}^{2}\right)  $ admits
a soliton metric which is isometric to $\mathbb{E}^{3}$. The singularity is
Type III and the Type III limit of compact manifolds modeled on $\widetilde
{\operatorname{Isom}}\left(  \mathbb{E}^{2}\right)  $ consist of compact,
collapsed quotients of $\mathbb{E}^{3}.$
\end{itemize}
\end{theorem}

\end{document}